\documentclass{amsart}
\usepackage{graphicx}
\usepackage{amssymb}
\usepackage{epstopdf}
\usepackage{nicefrac}
\usepackage{enumerate}
\usepackage{hyperref}
\usepackage{float}
\usepackage{color}
\usepackage{pst-all}
\usepackage{tabularx}

\newtheorem{thm}{THEOREM}[section]

\newtheorem{cor}[thm]{COROLLARY}

\newtheorem{defn}[thm]{DEFINITION}
\newtheorem{ex}[thm]{EXAMPLE}

\newtheorem{prop}[thm]{PROPOSITION}

\newtheorem{problem}[thm]{PROBLEM}
\newtheorem{remark}[thm]{REMARK}

\newcommand{\ds}{\displaystyle}

\newcommand{\cN}{{\mathcal N}}

\newcommand{\cD}{{\mathcal D}}

\newcommand{\cG}{{\mathcal G}}
\newcommand{\cH}{{\mathcal H}}
\newcommand{\cI}{{\mathcal I}}

\newcommand{\cO}{{\mathcal O}}
\newcommand{\cP}{{\mathcal P}}
\newcommand{\CO}{{\rm CO}} 
\newcommand{\cS}{{\mathcal S}}

\newcommand{\cU}{{\mathcal U}}

\newcommand{\diam}{{\rm diam}} 
 
\newcommand{\dX}{d_{\fX}} 
\newcommand{\e}{{\varepsilon}} 

\newcommand{\fD}{{\mathfrak{D}}}
\newcommand{\fH}{{\mathfrak{H}}}
\newcommand{\fG}{{\mathfrak{G}}}
\newcommand{\fK}{{\mathfrak{K}}}
\newcommand{\Fix}{{\rm Fix}} 
\newcommand{\fL}{{\mathfrak{L}}}
\newcommand{\fM}{{\mathfrak{M}}}
\newcommand{\fN}{{\mathfrak{N}}}

\newcommand{\fX}{{\mathfrak{X}}}

\newcommand{\fU}{{\mathfrak{U}}}
\newcommand{\G}{\Gamma}
\newcommand{\F}{{\mathcal F}}

\newcommand{\Homeo}{{\rm Homeo}}

\newcommand{\Iso}{{\rm Iso}} 

\newcommand{\mR}{{\mathbb R}}
\newcommand{\mS}{{\mathbb S}}
\newcommand{\mT}{{\mathbb T}}
\newcommand{\mZ}{{\mathbb Z}}

\newcommand{\Perm}{{\rm Perm}} 
\newcommand{\vp}{{\varphi}}
\newcommand{\whe}{\widehat{e}}
\newcommand{\whg}{\widehat{g}}

\newcommand{\whq}{\widehat{q}}

\newcommand{\whGamma}{\widehat{\Gamma}}

\newcommand{\whC}{\widehat{C}}
\newcommand{\whG}{\widehat{G}}
\newcommand{\whU}{\widehat{U}}
\newcommand{\whmZ}{\widehat{\mZ}}
\newcommand{\whPhi}{\widehat{\Phi}}

\newcommand{\whrho}{{\widehat{\rho}}}
\newcommand{\whtau}{{\widehat{\tau}}}
\newcommand{\whcH}{{\widehat{\cH}}}

\newcommand{\ovq}{\overline{q}}

\newcommand{\oa}{\overline{a}}
\newcommand{\ob}{\overline{b}}
\newcommand{\oc}{\overline{c}}

\newcommand{\ox}{\overline{x}}
\newcommand{\oy}{\overline{y}}
\newcommand{\oz}{\overline{z}}

\newcommand{\wha}{\widehat{a}}
\newcommand{\whb}{\widehat{b}}
\newcommand{\whc}{\widehat{c}}

\newcommand{\wtM}{{\widetilde{M}}}

 \newcommand\myeq{\mathrel{\stackrel{\makebox[0pt]{\mbox{\normalfont\tiny def}}}{=}}}

\newcommand\mor{\mathrel{\stackrel{\makebox[0pt]{\mbox{\normalfont\tiny a}}}{\sim}}}

\begin{document}

\title{The prime spectrum of solenoidal manifolds}

\author{Steven Hurder}
\address{Steven Hurder, Department of Mathematics, University of Illinois at Chicago, 322 SEO (m/c 249), 851 S. Morgan Street, Chicago, IL 60607-7045}
\email{hurder@uic.edu}

\author{Olga Lukina}
\address{Olga Lukina, Faculty of Mathematics, University of Vienna, Oskar-Morgenstern-Platz 1, 1090 Vienna, Austria}
\email{olga.lukina@univie.ac.at}

\thanks{Version date: March 11, 2021}


\thanks{OL is supported by FWF Project P31950-N35}


  \begin{abstract}
  A solenoidal manifold is the inverse limit space of a tower of proper coverings of a compact manifold. In this work, we introduce new   invariants for solenoidal manifolds, their asymptotic Steinitz orders and their prime spectra, and show they are invariants of the homeomorphism type. These invariants are formulated in terms of the monodromy Cantor action associated to a solenoidal manifold. To this end,   we continue our study of invariants for minimal equicontinuous Cantor actions.  We introduce the three types of prime spectra associated to such actions, and study their invariance properties under return equivalence.
  As an application, we show that a nilpotent Cantor action with finite prime spectrum must be stable. Examples of stable actions of the integer Heisenberg group are given with arbitrary prime spectrum. We also give the first examples of nilpotent Cantor actions  which  are wild, and not stable.
  \end{abstract}

\maketitle

 {\small  
   \tableofcontents
 }
 \vfill
 \eject
  
\section{Introduction}\label{sec-intro}

A 1-dimensional solenoid  is   the inverse limit space of a sequence of covering maps,  
  \begin{equation}\label{eq-presentation1}
\cS(\vec{m}) \myeq \underleftarrow{\lim} ~ \{\, q_{\ell}:\mS^1 \to \mS^1 \mid  \ell \geq 1\} 
\end{equation}
where $q_{\ell}$ is  a covering map  of the circle $\mS^1$   of degree  $m_{\ell} > 1$. Here,  $\vec{m} = (m_1, m_2,\dots)$  denotes a  sequence of   integers with each $m_i \geq  2$.  These continua (compact metric spaces) were introduced by van Danzig \cite{vanDantzig1930} and Vietoris \cite{Vietoris1927}, and appear  in many  areas of mathematics.

Associated to $\vec{m}$ is a \emph{supernatural} number, or \emph{Steinitz} number,  $\Pi[\vec{m}]$, 
which is the formal product of the integers $\{m_i \mid i \geq 1\}$.
 Chapter 2 of   Wilson  \cite{Wilson1998}, or Chapter~2.3 of Ribes and Zilesskii \cite{RZ2000},  give  a basic discussion of the  arithmetic of supernatural numbers. In particular, a Steinitz number 
   can be rewritten as the formal product of its prime factors, 
   \begin{equation}\label{eq-steinitzorder}
\Pi = \Pi[\vec{m}] = m_1 \cdot m_2 \cdots m_i \cdots  = \prod_{p \in \pi} \ p^{n(p)} \quad , \quad 0 \leq n(p) \leq \infty \ ,
\end{equation}
where 
 $\pi = \{2,3,5, \ldots\}$ is the set of distinct prime numbers. The non-negative integers $n(p)$ can be thought of as the ``coordinates'' of $\Pi$ along the ``axes'' given by the primes in $\pi$.

The Steinitz number $\Pi[\vec{m}]$ is called the \emph{Steinitz order} of the inverse limit $\cS(\vec{m})$.
  The following equivalence relation  appears naturally in the applications of Steinitz numbers to dynamical systems.

 \begin{defn}\label{def-asymptequiv}
Given $\vec{m} = \{m_i \mid i \geq 1 ~ , ~ m_i > 1 \}$ and $\vec{n} = \{n_i \mid i \geq 1 ~ , ~ n_i > 1\}$ sequences of integers, we say that the   Steinitz numbers $\Pi[\vec{m}]$ and $\Pi[\vec{n}]$ are \emph{asymptotically equivalent}, and we write   $\Pi[\vec{m}] \mor \Pi[\vec{n}]$, if there exist    integers $1 \leq m_0 < \infty$ and $1 \leq n_0 < \infty$ such that $n_0 \cdot \Pi[\vec{m}] = m_0 \cdot \Pi[\vec{n}]$.  The asymptotic equivalence class of $\Pi[\vec{m}]$ is denoted by  $\Pi_a[\vec{m}]$.
   \end{defn}
 Definition \ref{def-asymptequiv} says that two representatives of the same asymptotic equivalence class $\Pi_a[\vec{m}]$ differ by a finite number of prime factors with finite coordinates.
   
    Bing observed in \cite{Bing1960} that for $1$-dimensional solenoids $\cS(\vec{m})$ and $\cS(\vec{n})$, if  $\Pi[\vec{m}] \mor \Pi[\vec{n}]$ then the solenoids are homeomorphic.  McCord showed in \cite[Section~2]{McCord1965}   the converse, that if $\cS(\vec{m})$ and $\cS(\vec{n})$ are homeomorphic spaces, then  $\Pi[\vec{m}] \mor \Pi[\vec{n}]$.   Aarts and Fokkink  gave  in  \cite{AartsFokkink1991}  an alternate proof of this.  Thus we have:
    
\begin{thm}\cite{AartsFokkink1991,Bing1960}\label{thm-onedimSol}
Solenoids $\cS(\vec{m}) $ and $\cS(\vec{n}) $ are homeomorphic  if and only if ~ $\Pi[\vec{m}] \mor \Pi[\vec{n}]$. 
\end{thm}

The results in this paper  were motivated in part by the question, to what extent does Theorem~\ref{thm-onedimSol} generalize to higher dimensional solenoidal manifolds?  

A sequence of \emph{proper finite covering} maps 
$\cP = \{\, q_{\ell} \colon  M_{\ell} \to M_{\ell -1} \mid  \ell \geq 1\}$, where each $M_{\ell}$ is a compact connected manifold without boundary of dimension $n \geq 1$, is called a \emph{presentation} in \cite{DHL2016c}. The inverse limit     
\begin{equation}\label{eq-presentationinvlim}
\cS_{\cP} \equiv \lim_{\longleftarrow} ~ \{ q_{\ell } \colon M_{\ell } \to M_{\ell -1}\} ~ \subset \prod_{\ell \geq 0} ~ M_{\ell} ~  
\end{equation}
is   the  \emph{solenoidal manifold}   associated to $\cP$. The set $\cS_{\cP}$ is given  the relative  topology, induced from the product topology, so that $\cS_{\cP}$ is  compact and connected.  
By the definition of the inverse limit, for a sequence $\{x_{\ell} \in M_{\ell} \mid \ell \geq 0\}$, we have 
\begin{equation}\label{eq-presentationinvlim2}
x = (x_0, x_1, \ldots ) \in \cS_{\cP}   ~ \Longleftrightarrow  ~ q_{\ell}(x_{\ell}) =  x_{\ell-1} ~ {\rm for ~ all} ~ \ell \geq 1 ~. 
\end{equation}
For each $\ell \geq 0$, there  is a    fibration $\whq_{\ell} \colon \cS_{\cP} \to M_{\ell}$,  given by projection onto the $\ell$-th factor in \eqref{eq-presentationinvlim}, so $\whq_{\ell}(x) = x_{\ell}$.   We also make use of the covering maps denoted by
$\ovq_{\ell} = q_{\ell} \circ q_{\ell -1} \circ \cdots \circ q_1 \colon M_{\ell} \to M_0$. Note that $\whq_0 = \ovq_{\ell} \circ \whq_{\ell}$.
 
Solenoidal manifolds, as a special class of continua,   were first studied by McCord in \cite{McCord1965}, who showed that 
the continuum $\cS_{\cP}$ is a foliated space with foliation $\F_{\cP}$, in the sense of \cite{MS2006}, where the leaves of $\F_{\cP}$ are coverings of the base manifold $M_0$ via the projection map $\whq_0 \colon \cS_{\cP} \to M_0$ restricted to the path-connected components of  $\cS_{\cP}$. Solenoidal manifolds are     \emph{matchbox manifolds} of dimension $n$ in the terminology of \cite{ClarkHurder2013}, and the terminology    ``solenoidal manifolds'' was introduced by  Sullivan \cite{Sullivan2014}.     
The Heisenberg $H_3(\mR)$-odometers studied by  Danilenko and Lema\'{n}czyk in \cite{DL2016} are all solenoidal manifolds, equipped with the leafwise   action of $H_3(\mR)$.

The motivation for McCord's work in \cite{McCord1965} was the question of whether a solenoidal space must be a homogeneous continuum? That is, when does the group of self-homeomorphisms act transitively on the space? 
This is a particular case of the  more general problem      to study the space of homeomorphisms between solenoidal manifolds,  and their invariants up to homeomorphism. This problem has been studied especially in the  works   \cite{AartsFokkink1991,CHL2019,HL2018a}.
In this work, we continue this study by associating a \emph{prime spectrum} to a solenoidal space, and studying its invariance properties.

Given a presentation $\cP$,   define the truncated presentation 
$\cP_m = \{\, q_{\ell} \colon  M_{\ell} \to M_{\ell -1} \mid  \ell > m\}$, then it is a formality that the solenoidal manifolds $\cS_{\cP}$ and $\cS_{\cP_m}$ are homeomorphic. Thus, homeomorphism invariants for solenoidal manifolds   have an ``asymptotic'' character in terms of its presentation.

For a presentation $\cP$ as in \eqref{eq-presentationinvlim}, let $m_{\ell} > 1$ denote the degree of the covering map 
$q_{\ell } \colon M_{\ell } \to M_{\ell -1}$. The   product $m_{1} \cdots m_{\ell}$ equals the degree of the covering map $\ovq_{\ell} \colon M_{\ell} \to M_0$.
\begin{defn}\label{def-steinitzpres}
   The \emph{Steinitz order} of a presentation $\cP$ is the Steinitz number
   \begin{equation}\label{eq-highersteinitzorder}
\Pi[\cP] = LCM \{ m_1   m_2 \cdots m_{\ell} \mid  \ell > 0\} \ ,
\end{equation}
where LCM denotes the least common multiple of the collection of integers. The \emph{asymptotic Steinitz order} of $\cP$ is the class $\Pi_a[\cP]$ associated to $\Pi[\cP]$.
\end{defn}

That is, the Steinitz order of a presentation $\cP$ counts the number of appearances of distinct primes in the degrees of the covering maps $\ovq_{\ell} \colon M_{\ell} \to M_0$ for $\ell \geq 1$. Here $LCM$ should be understood in terms of Steinitz numbers, see Example \ref{ex-computeLCM} for more explanation.

Our first result is a direct generalization of one of the implications of Theorem~\ref{thm-onedimSol}.

\begin{thm}\label{thm-main1}
Let $\cS_{\cP}$ be a   solenoidal manifold  with     presentation $\cP$. Then the asymptotic order $\Pi_a[\cP]$ depends only on the homeomorphism type of $\cS_{\cP}$, and so defines the \emph{asymptotic Steinitz order}  of $\cS_{\cP}$  denoted by $\Pi_a[\cS_{\cP}]$.
\end{thm}
 
 Note that  McCord's proof in \cite[Section~2]{McCord1965} for $1$-dimensional solenoids uses Pontrjagin Duality, and his technique of proof     is only applicable for the case when the fundamental group of $M_0$ is abelian.
 
One cannot expect a converse to the conclusion of Theorem~\ref{thm-main1}  as in Theorem~\ref{thm-onedimSol}. 
For example, if $M_0 = \mT^n$ is the $n$-torus with $n > 1$,  Example~\ref{ex-toral}  constructs solenoidal manifolds over $\mT^n$ which have equal asymptotic orders,  but are not homeomorphic. Examples~\ref{ex-stable} and \ref{ex-wild} construct isospectral nilpotent   Cantor actions whose suspension solenoids   are not homeomorphic.

The proof of Theorem~\ref{thm-main1} is  based on the study of   the monodromy actions of solenoidal manifolds, and the fact that a homeomorphism between   solenoidal manifolds induces a return equivalence  between their global monodromy Cantor actions, as discussed in Section~\ref{subsec-morita}.  The Steinitz order invariants for  minimal equicontinuous  Cantor actions studied in this work are of independent interest, and will   be described next.

We say that $(\fX, \G, \Phi)$ is a \emph{ Cantor action} if $\G$ is a countable group, $\fX$ is a   Cantor   space, and $\Phi \colon \G \times \fX \to \fX$ is a minimal  action.  The action  $(\fX,\G,\Phi)$ is \emph{equicontinuous} with respect to a metric $\dX$ on $\fX$, if for all $\e >0$ there exists $\delta > 0$, such that for all $x , y \in \fX$   with 
 $\ds  \dX(x,y) < \delta$ and  all  $\gamma \in \G$,    we have  $\dX(\gamma   x, \gamma   y) < \e$.
This property   is independent of the choice of the metric   on $\fX$.

Let  $\Phi(\G) \subset \Homeo(\fX)$ denote the image subgroup for an action $(\fX, \G, \Phi)$.
When the action is equicontinuous,     the  closure $\overline{\Phi(\G)} \subset \Homeo(\fX)$ in the \emph{uniform topology of maps} is a separable profinite group.  We adopt the notation $\fG(\Phi) \equiv \overline{\Phi(\G)}$. More generally, we typically use   letters in fraktur font   to denote profinite objects.     Let $\whPhi \colon \fG(\Phi) \times \fX \to \fX$ denote the induced   action of $\fG(\Phi)$ on $\fX$, which is transitive as the action   $(\fX, \G, \Phi)$ is   minimal. For $\whg \in \fG(\Phi)$, we   write its action on $\fX$ by $\whg \, x = \whPhi(\whg)(x)$.
Given $x \in \fX$,   introduce the isotropy group  at $x$,  
\begin{align}\label{iso-defn2}
 \fD(\Phi, x) = \{ \whg  \in \fG(\Phi) \mid \whg \, x = x\} \subset \Homeo(\fX) \ ,
\end{align}
which is a closed subgroup of $\fG(\Phi)$, and    thus is either finite, or is an infinite profinite group.    
 As the action $\whPhi \colon \fG(\Phi) \times \fX \to \fX$ is transitive,  the conjugacy class of $\fD(\Phi,x)$ in $\fG(\Phi)$ is independent of the choice of   $x$.     The group $\fD(\Phi,x)$ is called the \emph{discriminant} of the action $(\fX,\G,\Phi)$ in the authors works \cite{DHL2016c,HL2018a,HL2018b}, and is called a \emph{parabolic} subgroup  (of the profinite completion of a countable group)  in the works by Bartholdi and Grigorchuk \cite{BartholdiGrigorchuk2000,BartholdiGrigorchuk2002}.

 The \emph{Steinitz order} $\Pi[\fG]$ of a   profinite group $\fG$ is  a supernatural number associated to a presentation of $\fG$ as an inverse limit of finite groups (see  Definition~\ref{def-steinitzprofinite}, or   \cite[Chapter 2]{Wilson1998} or \cite[Chapter~2.3]{RZ2000}).
 The Steinitz order has been used in the study of the analytic representations of profinite groups associated to groups acting on rooted  trees, for example  in the work \cite{Kionke2019}. Parabolic subgroups of countable groups, acting on rooted trees, play an important role in the study of analytic representations of such groups, see for instance \cite{BartholdiGrigorchuk2000,BartholdiGrigorchuk2002}, and the importance of developing a similar theory for representations of profinite groups was pointed out in \cite{BartholdiGrigorchuk2002}.

Recall that for  a profinite group $\fG$, an open subgroup $\fU \subset \fG$ has finite index  \cite[Lemma 2.1.2]{RZ2000}.  Given a collection of finite positive integers $S = \{n_i \mid i \in \cI\}$, let  $LCM(S)$ denote the least common multiple of the collection, in the sense of Steinitz numbers.
\begin{defn}\label{def-steinitzorderaction}
Let $(\fX, \G, \Phi)$ be a minimal equicontinuous Cantor action, with choice of a basepoint $x \in \fX$. The \emph{Steinitz orders} of the action  are defined as follows:
 \begin{enumerate}
\item $\Pi[\fG(\Phi)] =  LCM \{\# \ \fG(\Phi)/\fN  \mid \fN \subset \fG(\Phi)~ \text{open normal subgroup}\}$, 
\item $\Pi[\fD(\Phi)] =  LCM \{\# \ \fD(\Phi,x)/(\fN \cap \fD(\Phi,x)) \mid \fN \subset \fG(\Phi)~ \text{open normal subgroup}\}$, 
\item $\Pi[\fG(\Phi) : \fD(\Phi)] =  LCM \{\# \ \fG(\Phi)/(\fN \cdot \fD(\Phi,x))  \mid  \fN \subset \fG(\Phi)~ \text{open normal subgroup}\}$.
\end{enumerate}
\end{defn}

The next result shows that these Steinitz orders are invariants of the isomorphism class of the action, for the notion of isomorphism or conjugacy as   given in Definition~\ref{def-isomorphism}.
 
 \begin{thm}\label{thm-isoinvariance}
Let  $(\fX, \G, \Phi)$  be a minimal equicontinuous Cantor action. Then the Steinitz orders for the action are independent of the choice of a basepoint $x \in \fX$.
Moreover, these orders depend only on the isomorphism class of the action, and satisfy the Lagrange identity
 \begin{equation}\label{eq-productorders}
\Pi[\fG(\Phi)]  = \Pi[\fG(\Phi) : \fD(\Phi)] \cdot \Pi[\fD(\Phi)]  \ ,
\end{equation}
where the multiplication is taken in the sense of supernatural numbers. 
\end{thm}

For example, if $\Phi \colon \mZ \times \fX \to \fX$ is a minimal equicontinuous action of the free abelian group $\G = \mZ$, which is the monodromy of a solenoid $\cS(\vec{m}) $ as  defined by \eqref{eq-presentation1}, then the Steinitz order of the closure of the action is given by $\Pi[\fG(\Phi)] = \Pi[\vec{m}]$. As the group $\G = \mZ$ is abelian, the discriminant subgroup $\fD(\Phi)$ is trivial, so $\Pi[\fD(\Phi)]$ is trivial, and $\Pi[\fG(\Phi) : \fD(\Phi)] = \Pi[\fG(\Phi)]$. On the other hand, there are   Cantor  actions of the   Heisenberg group with $\fD(\Phi)$ a  Cantor group, and their  Steinitz orders $[\fD(\Phi)]$ distinguish an uncountable number of such actions. (See the examples  in Section~\ref{subsec-heisenberg}.)

Isomorphism  is the strongest notion of equivalence for Cantor actions.   \emph{Return equivalence}, as given in Definition~\ref{def-return},  is a form of ``virtual   isomorphism'' for minimal equicontinuous Cantor actions, and   is natural when considering   Cantor systems arising from geometric constructions, as    in    \cite{HL2018a,HL2018b,HL2019a}.

 \begin{thm}\label{thm-returnequivorder}
 Let $(\fX, \G, \Phi)$ be a minimal equicontinuous   Cantor action. Then the relative asymptotic Steinitz order 
 $\Pi_a[\fG(\Phi) : \fD(\Phi)]$ is an invariant of   its return equivalence class.
\end{thm}
It is shown in Section~\ref{subsec-solenoids} that the Steinitz number $\Pi[\cP]$ of a presentation  in Theorem~\ref{thm-main1} equals the relative Steinitz order 
 $\Pi[\fG(\Phi) : \fD(\Phi)]$ for the monodromy action of the solenoid $\cS_{\cP}$, so that Theorem~\ref{thm-main1}  follows from Theorem~\ref{thm-returnequivorder} and the results of Sections~\ref{subsec-equivalence} and \ref{subsec-morita}.
 
The behavior under return equivalence of actions of the other two Steinitz orders $\Pi[\fG(\Phi)]$ and $\Pi[\fD(\Phi)]$ in Definition~\ref{def-steinitzorderaction}   is more subtle.   In particular, the constructions in  Example~\ref{ex-almosttoral} show that their asymptotic classes need not be invariant under return equivalence.   

  \begin{defn}\label{def-primespectrum}
Let  $\pi = \{2,3,5, \ldots\}$ denote the set of   primes. Given  $\Pi = \prod_{p \in \pi} \ p^{n(p)}$, define:
\begin{eqnarray*}
\pi(\Pi) ~ & = & ~  \{ p \in \pi \mid 0 < n(p)   \}   \ , ~ \emph{the prime spectrum of} \ \Pi ,  \label{eq-primespectrum}\\
\pi_f(\Pi) ~ & = & ~  \{ p \in \pi \mid 0 < n(p) < \infty \} \ , ~  \emph{the finite prime spectrum of} \ \Pi ,  \label{eq-finiteprimespectrum} \\
\pi_{\infty}(\Pi) ~ & = & ~  \{ p \in \pi \mid  n(p) = \infty \} \ , ~  \emph{the infinite prime spectrum of} \ \Pi \ .\label{eq-infiniteprimespectrum}
\end{eqnarray*}
   \end{defn}
  Note that if $\Pi \mor \Pi'$, then $\pi_{\infty}(\Pi) = \pi_{\infty}(\Pi')$. The property that $\pi_f(\Pi)$ is a \emph{infinite} set is also preserved by asymptotic equivalence of Steinitz numbers.

 A profinite group  $\fG$ is said to have \emph{finite prime spectrum} if ~  $\pi(\Pi(\fG))$ is a finite set of primes.  
If $\pi(\Pi(\fG)) = \{p\}$, then  $\fG$ is said to be a \emph{pro-$p$ group}, for which there is an extensive literature  \cite{DdSMS1999,dSSS2000}.
The property that $\Pi(\fG)$ has finite prime spectrum is preserved by   asymptotic equivalence.
   
  \begin{thm}\label{thm-returnequivspectra}
Let $(\fX, \G, \Phi)$ be a minimal equicontinuous   Cantor action. Then the  infinite prime spectra of the Steinitz orders 
 $\Pi[\fG(\Phi)]$,  $\Pi[\fD(\Phi)]$  and $\Pi_a[\fG(\Phi) : \fD(\Phi)]$  depend only on  the return equivalence class of the action.
The same holds for the property  that the finite prime spectrum of each of these Steinitz orders   is an  infinite set.
 \end{thm}

This result suggests a natural question:
 \begin{problem}\label{prob-nilpotentorder}
How do the dynamical properties of a minimal equicontinuous  Cantor action  $(\fX, \G, \Phi)$  depend on the asymptotic Steinitz orders associated to the action?
 \end{problem}

A basic dynamical property of a minimal equicontinuous  Cantor action $(\fX, \G, \Phi)$ is its degree of ``regularity'', as discussed in Section \ref{subsec-lqa}.  The action is \emph{topologically free} if the set of all fixed points for the elements of the action is a meagre set (see Definition~\ref{def-topfree}.) The \emph{local  quasi-analytic} property of an action, as in Definition~\ref{def-LQA}, is a local (generalized) version of the topologically free property, and does not require that the acting group $\G$ be countable, so applies for profinite group actions in particular.   We then have the following notion:
   \begin{defn}  \label{def-stable}  
 An equicontinuous  Cantor action  $(\fX, \G, \Phi)$ is said to be \emph{stable} if the  induced  profinite action  $\whPhi \colon \fG(\Phi)   \times \fX \to \fX$ is locally quasi-analytic. The action  is said to be \emph{wild} otherwise. 
\end{defn}
A stable Cantor action satisfies local rigidity, as discussed in the works \cite{CortezMedynets2016,HL2018b,HL2020a,Li2018}. On the other hand, 
there are many examples of wild Cantor actions. The actions of  weakly branch groups on the boundaries of their associated trees are always wild  \cite{BGS2012,GL2019}.  The work \cite{ALBLLN2020} gives the construction of wild Cantor actions exhibiting a variety of characteristic properties, using algebraic methods. 

 In this work, we a  partial solution  to Problem~\ref{prob-nilpotentorder}.
A \emph{nilpotent   Cantor action}  is a  minimal equicontinuous Cantor action  $(\fX,\G,\Phi)$,  where $\G$ contains a finitely-generated nilpotent subgroup $\G_0 \subset \G$ of finite index.   The authors showed in \cite[Theorem~4.1]{HL2020a}     that a nilpotent Cantor action is always locally quasi-analytic. Moreover, it was shown in \cite[Theorem~1.1]{HL2020a} that if the actions are both effective, then the property of being a nilpotent Cantor action is preserved by return equivalence, and thus also by continuous orbit equivalence of actions.

 \begin{thm}\label{thm-nilstable}
Let $(\fX,\G,\Phi)$ be a nilpotent Cantor action, with discriminant   $\fD(\Phi) \subset \fG(\Phi)$. 
If the prime spectrum $\pi(\Pi(\fD(\Phi)))$ is   finite,  then the action is stable. In particular, if the prime spectrum $\pi(\Pi[\fG(\Phi)])$ is   finite,  then the action is stable.
\end{thm}
 The proof of Theorem \ref{thm-nilstable}    yields the following corollary.
 The \emph{multiplicity} of a prime $p$ in a Steinitz number $\Pi$ is the value of $n(p)$ in the formula \eqref{eq-steinitzorder}. 
   \begin{cor}\label{cor-niltopfree}
Let $(\fX,\G,\Phi)$ be a nilpotent Cantor action. If the Steinitz order $\Pi(\fG(\Phi))$ has prime multiplicities at most $2$, for all but a finite set of primes, then the action is stable.
\end{cor}
The wild actions   in Example~\ref{ex-wild} have finite multiplicities at least $3$ for an infinite set of primes.

 The converse of Theorem \ref{thm-nilstable} need not hold, indeed, it is possible to construct actions of abelian groups with infinite prime spectrum which are necessarily stable, see Example \ref{ex-toral}, and also stable actions of nilpotent groups with infinite prime spectrum, see Example \ref{ex-stable}. The relation of the finite prime spectrum with the stability of an action depends on the \emph{Noetherian} property of its profinite completion, as explained in Section \ref{subsec-Noetheriandynamics}. 

The celebrated Grigorchuk group (see \cite{BGS2012,Grigorchuk2011} for example)  is a $p$-group for $p=2$, and its action on the boundary of the $2$-adic tree is minimal and equicontinuous, and moreover is a wild action.   Thus,    Theorem~\ref{thm-nilstable}  cannot be   generalized to Cantor actions of arbitrary finitely generated groups. 

 The authors asked in the works \cite{HL2018b,HL2020a} whether a locally quasi-analytic nilpotent Cantor action $(\fX,\G,\Phi)$ can be wild,  more precisely, do there exist actions $(\fX,\G,\Phi)$ such that the action of $\G$ on $\fX$ is locally quasi-analytic, while the action of the completion $\fG(\Phi)$ on $\fX$ is not locally quasi-analytic? Using the constructions in Example~\ref{ex-wild}, our final result gives an   answer to this question, noting that a topologically-free Cantor action is locally quasi-analytic.

    \begin{thm}\label{thm-lqaWILD}
  There exists an uncountable number of topologically-free Cantor actions $(\fX,\G,\Phi)$ of the Heisenberg group $\G$, distinct up to return equivalence, that  are wild. 
  \end{thm}

 Section~\ref{sec-basics} recalls some basic facts about Cantor actions as required for this work. 
 
 Section~\ref{sec-supernatural} develops in more detail the properties of Steinitz orders for Cantor actions. This yield the proofs of 
  Theorems~\ref{thm-isoinvariance}, \ref{thm-returnequivorder} and \ref{thm-returnequivspectra}. Then in 
 Section~\ref{subsec-algmodel} we recall  the construction of the group chain model for a minimal equicontinuous Cantor action, and the results of Section~\ref{subsec-steinitzalgmodel} show that their Steinitz orders can be calculated using these group chains. This is used to deduce  the proof of Theorem~\ref{thm-main1} from Theorem~\ref{thm-returnequivorder} in Section~\ref{subsec-solenoids}.
 
Section~\ref{sec-nilpotent} considers the special case of nilpotent Cantor actions, and gives an application of the prime spectrum to this class of actions. 

 An essential part of the abstract study of minimal equicontinuous Cantor actions is to have explicit examples of the properties being studied and characterized. This we provide in Section~\ref{sec-examples}.
 
 Example~\ref{ex-toral} gives the most basic construction of   actions with prescribed prime spectrum for $\fG(\Phi)$. The $\mZ^n$-actions constructed in   show that for $n \geq 2$,    the prime spectrum does not contain sufficient information about the action to distinguish the actions up to return equivalence.

Example~\ref{ex-trivial} recalls the construction from \cite{HLvL2020} of a ``balanced"   self-embedding of the integer Heisenberg group into itself, which has the property that the discriminant group $\fD(\Phi)$ of the action is trivial, but the maps in the inverse limit formula for $\fD(\Phi)$ in  \eqref{eq-discformula} are not surjective.

  Example~\ref{ex-stable}   gives the construction of nilpotent Cantor actions of the integer Heisenberg group with arbitrary finite  or infinite prime spectrum, for  which the discriminant group $\fD(\Phi)$ is non-trivial and the action is stable.  
  Example~\ref{ex-wild}   gives  the constructions of   nilpotent Cantor actions for which the prime spectrum is any arbitrary infinite subset of the primes,   and   the action is   wild. These examples are then used to give the proof of Theorem~\ref{thm-lqaWILD}.

 \section{Cantor actions}\label{sec-basics}

We recall some of the basic 
 properties of     Cantor actions, as required for the proofs of the results in Section~\ref{sec-intro}.
More complete discussions of the properties of equicontinuous Cantor actions are given in     the text by Auslander \cite{Auslander1988}, the papers by Cortez and Petite  \cite{CortezPetite2008}, Cortez and Medynets  \cite{CortezMedynets2016},    and   the authors' works, in particular   \cite{DHL2016c} and   \cite[Section~3]{HL2019a}.

\subsection{Basic concepts}\label{subsec-basics}

Let  $(\fX,\G,\Phi)$   denote an action  $\Phi \colon \G \times \fX \to \fX$. We   write $g\cdot x$ for $\Phi(g)(x)$ when appropriate.
The orbit of  $x \in \fX$ is the subset $\cO(x) = \{g \cdot x \mid g \in \G\}$. 
The action is \emph{minimal} if  for all $x \in \fX$, its   orbit $\cO(x)$ is dense in $\fX$.

 \eject

   Let $N(\Phi) \subset \G$ denote the kernel of the action homomorphism $\Phi \colon \G \to \Homeo(\fX)$. The action is said to be \emph{effective} if $N(\Phi)$ is the trivial group. That is, the homomorphism $\Phi$ is faithful, and one also says that the action  is faithful.
    
 An action  $(\fX,\G,\Phi)$ is \emph{equicontinuous} with respect to a metric $\dX$ on $\fX$, if for all $\e >0$ there exists $\delta > 0$, such that for all $x , y \in \fX$ and $g \in \G$ we have  that 
 $\ds  \dX(x,y) < \delta$ implies   $\dX(g \cdot x, g \cdot y) < \e$.
The property of being equicontinuous    is independent of the choice of the metric   on $\fX$ which is   compatible with the topology of $\fX$.

Now assume that $\fX$ is a Cantor space. 
Let $\CO(\fX)$ denote the collection  of all clopen (closed and open) subsets of  $\fX$, which forms a basis for the topology of $\fX$. 
For $\phi \in \Homeo(\fX)$ and    $U \in \CO(\fX)$, the image $\phi(U) \in \CO(\fX)$.  
The following   result is folklore, and a proof is given in \cite[Proposition~3.1]{HL2018b}.
 \begin{prop}\label{prop-CO}
 For $\fX$ a Cantor space, a minimal   action   $\Phi \colon \G \times \fX \to \fX$  is  equicontinuous  if and only if  the $\G$-orbit of every $U \in \CO(\fX)$ is finite for the induced action $\Phi_* \colon \G \times \CO(\fX) \to \CO(\fX)$.
\end{prop}

We say that $U \subset \fX$  is \emph{adapted} to the action $(\fX,\G,\Phi)$ if $U$ is a   \emph{non-empty clopen} subset, and for any $g \in \G$, 
if $\Phi(g)(U) \cap U \ne \emptyset$ implies that  $\Phi(g)(U) = U$.   The proof of   \cite[Proposition~3.1]{HL2018b} shows that given  $x \in \fX$ and clopen set $x \in W$, there is an adapted clopen set $U$ with $x \in U \subset W$. 

For an adapted set $U$,   the set of ``return times'' to $U$, 
 \begin{equation}\label{eq-adapted}
\G_U = \left\{g \in \G \mid g \cdot U  \cap U \ne \emptyset  \right\}  
\end{equation}
is a subgroup of   $\G$, called the \emph{stabilizer} of $U$.      
  Then for $g, g' \in \G$ with $g \cdot U \cap g' \cdot U \ne \emptyset$ we have $g^{-1} \, g' \cdot U = U$, hence $g^{-1} \, g' \in \G_U$. Thus,  the  translates $\{ g \cdot U \mid g \in \G\}$ form a finite clopen partition of $\fX$, and are in 1-1 correspondence with the quotient space $X_U = \G/\G_U$. Then $\G$ acts by permutations of the finite set $X_U$ and so the stabilizer group $\G_U \subset G$ has finite index.  Note that this implies that if $V \subset U$ is a proper inclusion of adapted sets, then the inclusion $\G_V \subset \G_U$ is also proper.

\begin{defn}\label{def-adaptednbhds}
Let  $(\fX,\G,\Phi)$   be a minimal equicontinuous Cantor    action.
A properly descending chain of clopen sets $\cU = \{U_{\ell} \subset \fX  \mid \ell \geq 0\}$ is said to be an \emph{adapted neighborhood basis} at $x \in \fX$ for the action $\Phi$,   if
    $x \in U_{\ell +1} \subset U_{\ell}$  is a proper inclusion for all $ \ell \geq 0$, with     $\cap_{\ell > 0}  \ U_{\ell} = \{x\}$, and  each $U_{\ell}$ is adapted to the action $\Phi$.
\end{defn}
Given $x \in \fX$ and   $\e > 0$, Proposition~\ref{prop-CO} implies there exists an adapted clopen set $U \in \CO(\fX)$ with $x \in U$ and $\diam(U) < \e$.  Thus, one can choose a descending chain $\cU$ of adapted sets in $\CO(\fX)$ whose intersection is $x$, from which the following result follows:

\begin{prop}\label{prop-adpatedchain}
Let  $(\fX,\G,\Phi)$   be a minimal equicontinuous Cantor    action. Given $x \in \fX$, there exists an adapted neighborhood basis $\cU$ at $x$ for the action $\Phi$.
 \end{prop}

\subsection{Equivalence of Cantor actions}\label{subsec-equivalence}

We next recall the   notions of equivalence of  Cantor actions which we use in this work. 
The first and strongest   is that  of 
  {isomorphism} of Cantor actions, which   is a   generalization  of the usual notion of conjugacy of topological actions. For $\G = \mZ$, isomorphism corresponds to the notion of ``flip conjugacy'' introduced    in the work of Boyle and Tomiyama \cite{BoyleTomiyama1998}.
The definition below agrees with the usage in the papers     \cite{CortezMedynets2016,HL2018b,Li2018}.

 \begin{defn} \label{def-isomorphism}
Cantor actions $(\fX_1, \G_1, \Phi_1)$ and $(\fX_2, \G_2, \Phi_2)$ are said to be \emph{isomorphic}  if there is a homeomorphism $h \colon \fX_1 \to \fX_2$ and group isomorphism $\Theta \colon \G_1 \to \G_2$ so that 
\begin{equation}\label{eq-isomorphism}
\Phi_1(g) = h^{-1} \circ \Phi_2(\Theta(g)) \circ h   \in   \Homeo(\fX_1) \   \textrm{for  all} \ g \in \G_1 \ .
\end{equation}
 \end{defn}

The notion of \emph{return equivalence} for Cantor actions is  weaker than the notion of isomorphism, and is natural when considering the Cantor systems defined by the holonomy actions for solenoidal manifolds, as considered in   the works  \cite{HL2018a,HL2018b,HL2019a}.

For a minimal equicontinuous Cantor action $(\fX, \G, \Phi)$ and   an adapted set $U \subset \fX$, by a small abuse of  notation, we use $\Phi_U$ to denote both the restricted action $\Phi_U \colon \G_U \times U \to U$ and the induced quotient action $\Phi_U \colon H_U \times U \to U$ for $H_U = \Phi(G_U) \subset \Homeo(U)$. Then $(U, H_U, \Phi_U)$ is called the \emph{holonomy action} for $\Phi$, in analogy with the case where $U$ is a transversal to a solenoidal manifold, and 
  $H_U$ is the holonomy group for this transversal.

   \begin{defn}\label{def-return}
Two minimal equicontinuous Cantor  actions $(\fX_1, \G_1, \Phi_1)$ and $(\fX_2, \G_2, \Phi_2)$ are  \emph{return equivalent} if there exists 
  an adapted set $U_1 \subset \fX_1$ for the action $\Phi_1$   and  
  an adapted set $U_2 \subset \fX_2$ for the action $\Phi_2$,
such that   the  restricted actions $(U_1, H_{1,U_1}, \Phi_{1,U_1})$ and $(U_2, H_{2,U_2}, \Phi_{2,U_2})$ are isomorphic.
\end{defn}
If the actions $\Phi_1$ and $\Phi_2$ are isomorphic in the sense of Definition~\ref{def-isomorphism}, then they are return equivalent with   $U_1 = \fX_1$ and $U_2 = \fX_2$. However, the notion of return equivalence is weaker even for this case, as the conjugacy is between the holonomy groups $H_{1,\fX_1}$ and $H_{2,\fX_2}$, and not the groups $\G_1$ and $\G_2$.

\subsection{Morita equivalence}\label{subsec-morita}

We next relate the notion of return equivalence  of Cantor actions with that of Morita equivalence of pseudogroups, as induced by a homeomorphism  between solenoidal manifolds. Let $h \colon \cS_{\cP} \to \cS_{\cP'}$ be a homeomorphism between solenoidal manifolds, defined by 
$$\cS_{\cP} \equiv \lim_{\longleftarrow} ~ \{ q_{\ell } \colon M_{\ell } \to M_{\ell -1}\} ~ \subset \prod_{\ell \geq 0} ~ M_{\ell} \quad  , \quad 
\cS_{\cP'} \equiv \lim_{\longleftarrow} ~ \{ q'_{\ell } \colon M'_{\ell } \to M'_{\ell -1}\} ~ \subset \prod_{\ell \geq 0} ~ M'_{\ell} ~ ,
$$
with foliations $\F_{\cP}$ and $\F_{\cP'}$ defined by the path-connected components of each space, respectively.

Let  $\whq_0 \colon \cS_{\cP} \to M_0$  and  $\whq_0' \colon \cS_{\cP'} \to M'_0$  be the corresponding projection maps.
Then for choices of basepoints $x \in \cS_{\cP}$ and $x' \in \cS_{\cP'}$,   the Cantor fibers $\fX = \whq_0^{-1}(\whq_0(x))$  and $\fX' = (\whq'_0)^{-1}(\whq'_0(x'))$ are complete transversals to the foliations $\F_{\cP}$ and $\F_{\cP'}$, respectively.
The homeomorphism $h$ cannot be assumed to be fiber-preserving; that is, to satisfy $h(\fX) = \fX'$. For example, the work \cite{CHL2019} studies the homeomorphisms between solenoidal manifolds induced by    lifts of   homeomorphisms between finite covering spaces  $\pi \colon \wtM_0 \to M_0$ and $\pi' \colon \wtM_0' \to M_0'$ in which case the map $h$ need not even be continuously deformable into a fiber-preserving map.

 Associated to the transversal $\fX$ for $\F_{\cP}$ is a pseudogroup $\cG$ modeled on $\fX$. The elements of  $\cG$ are local homeomorphisms between open subsets $U,V \subset \fX$ induced by the holonomy transport along the leaves of $\F_{\cP}$. The construction of these pseudogroups for smooth foliations is discussed by Haefliger in \cite{Haefliger1984,Haefliger2002a}, for example. The adaptation of these ideas to matchbox manifolds, where the transverse space is a Cantor set, is discussed in detail in the works \cite{ClarkHurder2013,CHL2019}. 
 
 Associated to a non-empty open subset $W \subset \fX$, we can form the restricted pseudogroup $\cG_W$ which consists of the elements of $\cG$ whose domain and range are contained in $W$. As the foliation $\F_{\cP}$ is minimal, that is, every leaf is dense in $\cS_{\cP}$,  the pseudogroups $\cG$ and $\cG_W$ are Morita equivalent in the sense of Haefliger in \cite{Haefliger1984}.
 The same remarks apply to the space $\cS_{\cP'}$ and so there is a restricted pseudogroup $\cG'_{W'}$ for the pseudogroup $\cG'$ modeled on $\fX'$ defined by the holonomy transport of $\F_{\cP'}$.
 
The homeomorphism $h \colon \cS_{\cP} \to \cS_{\cP'}$ is  necessarily leaf-preserving, and a basic fact is that there exists non-empty open sets $W \subset \fX$ and $W' \subset \fX'$ such that  the homeomorphism $h$ induces an isomorphism between the restricted pseudogroups $\cG_W$ and $\cG'_{W'}$. This is discussed in detail in \cite[Section~2.4]{HL2018a}. Moreover, as the holonomy action of $\cG$ on $\fX$ is equicontinuous, and likewise that for $\cG'$ on $\fX'$, the open sets $W$ and $W'$ can be chosen to be clopen. Moreover,  $\cG_W$ is the pseudogroup induced by a minimal equicontinuous  group action on $W$, and likewise for the action of $\cG'_{W'}$ on $W'$, so $h$ induces a return equivalence between these group actions in the sense of Definition~\ref{def-return}.  Then by the remarks in Section~\ref{subsec-solenoids}, the algebraic model Cantor actions for the monodromy actions of $\cS_{\cP}$ and $\cS_{\cP'}$ are return equivalent.

 \subsection{Regularity of Cantor actions}\label{subsec-lqa}
 We next recall some regularity properties of Cantor actions. These are used in   the proof of Theorem~\ref{thm-nilstable} and the analysis of the examples constructed in Section~\ref{sec-examples}.

An action  $(\fX,\G,\Phi)$ is   said to be  \emph{free} if for all $x \in \fX$ and $g \in \G$,   $g \cdot x = x$ implies that $g = e$,   the identity of the group. 
   The   notion of a \emph{topologically free} action is a generalization of free actions,    introduced by Boyle in his thesis \cite{Boyle1983}, and later  used  in the works by Boyle and Tomiyama \cite{BoyleTomiyama1998}   for the study of classification of general Cantor actions,     by   Renault \cite{Renault2008}   for the study of the $C^*$-algebras associated to Cantor actions, and by   Li \cite{Li2018} for proving    rigidity properties of Cantor actions. We recall this definition.

Let $\Fix(g) = \{x \in \fX \mid g \cdot x = x \}$, and define the \emph{isotropy set}
\begin{equation}\label{eq-isotropy}
 \Iso(\Phi) = \{ x \in \fX \mid \exists ~ g \in \G ~ , ~ g \ne id ~, ~g \cdot x = x    \} = \bigcup_{e \ne g \in \G} \ \Fix(g) \ . 
\end{equation}

  \begin{defn}\cite{BoyleTomiyama1998,Li2018,Renault2008} \label{def-topfree}
  $(\fX,\G,\Phi)$ is said to be \emph{topologically free}  if  $\Iso(\Phi) $ is meager in $\fX$. 
 \end{defn}
 
Note that if $\Iso(\Phi)$ is meager, then $\Iso(\Phi)$ has empty interior. That is, if there exists a non-identity element $g \in \G$ such that $\Fix(g)$ has interior, then the action is not topologically free.

The notion of a 
  \emph{quasi-analytic} action,  introduced in the  works  of    {\'A}lvarez L{\'o}pez,   Candel, and Moreira Galicia  \cite{ALC2009,ALM2016}, is an alternative formulation of the topologically free property which generalizes to group Cantor actions   where the acting group can be countable or profinite. 
  \begin{defn}\label{def-qa}
 An action $\Phi \colon H \times \fX \to \fX$, where 
    $H$ is a topological group and  $\fX$ a Cantor space,    is said to be \emph{quasi-analytic} if for each clopen set $U \subset \fX$  
  and $g \in H$ such that $\Phi(g)(U) = U$ and the restriction $\Phi(g) | U$ is the identity map on $U$, 
  then $\Phi(g)$ acts as the identity on  $\fX$.  
  \end{defn}

A topologically free action  is quasi-analytic.   Conversely,  the Baire Category Theorem implies that a quasi-analytic   effective action of a \emph{countable} group  is topologically free \cite[Section~3]{Renault2008}.

  A local formulation of the quasi-analytic property was introduced  in the works \cite{DHL2016c,HL2018a}, and has proved very useful for the study of the dynamical properties of Cantor actions. 
    \begin{defn}  \label{def-LQA}  
   An action $\Phi \colon H \times \fX \to \fX$, where 
    $H$ is a topological group and  $\fX$ a Cantor metric space with metric $\dX$,   is   \emph{locally quasi-analytic}  (or LQA) if there exists $\e > 0$ such that for any non-empty open set $U \subset \fX$ with $\diam (U) < \e$,  and  for any non-empty open subset $V \subset U$,  if the action of $g \in H$ satisfies $\Phi(g)(V) = V$ and the restriction $\Phi(g) | V$ is the identity map on $V$,    then $\Phi(g)$ acts as the identity on all of $U$.  
\end{defn}
 
 This reformulation of  the notion of   topologically free actions is the basis for the following notion.
 \begin{defn}\label{def-stable2}
 A minimal equicontinuous Cantor action  $(\fX, \G, \Phi)$ is said to be \emph{stable} if the action of its profinite closure $\fG(\Phi)$ on $\fX$ is locally quasi-analytic, and otherwise is  a \emph{wild} action.
 \end{defn}

Wild Cantor actions include the actions of  weakly  branch groups on their boundaries 
\cite{BartholdiGrigorchuk2000,BartholdiGrigorchuk2002,BGS2012,DudkoGrigorchuk2017,Grigorchuk2011,Nekrashevych2005,Nekrashevych2016},   actions of higher rank arithmetic lattices on quotients of their profinite completions \cite{HL2018a}, and various constructions of subgroups of wreath product groups acting on trees \cite{ALBLLN2020}.

 \section{Steinitz orders of Cantor actions}\label{sec-supernatural}

In this section, we  recall  the properties of the Steinitz orders of profinite groups from the texts \cite{RZ2000,Wilson1998}, then consider the invariance properties of the Steinitz orders associated to a minimal equicontinuous Cantor action. This yields    
   proofs of Theorems~\ref{thm-isoinvariance}, \ref{thm-returnequivorder} and \ref{thm-returnequivspectra}.
  We then recall the algebraic model for a minimal equicontinuous action, and derive the Steinitz orders of a Cantor action  in terms of this algebraic model. The algebraic models are used in the proof of Theorem~\ref{thm-main1} in Section~\ref{subsec-solenoids}, and for the constructions of examples in Section~\ref{sec-examples}.

\subsection{Abstract Steinitz orders}\label{subsec-sodefs} 
We begin with the definitions and basic properties of the Steinitz orders associated to profinite groups.
\begin{defn}\label{def-steinitzprofinite}
 Let $\fH \subset \fG$ be a closed subgroup  of the profinite group $\fG$. Then  
  \begin{equation}\label{eq-relativeorder}
\Pi[\fG : \fH] =   LCM \{\# \ \fG(\Phi)/(\fN \cdot \fH)  \mid \fN \subset \fG(\Phi)~ \text{clopen normal subgroup}\}  
\end{equation}
is the \emph{relative Steinitz order} of $\fH$ in $\fG$. The \emph{Steinitz order} of $\fG$ is 
$\Pi[\fG] = \Pi[\fG : \{\whe\}]$, where   $\{\whe\}$ is the identity subgroup. 
\end{defn}
  \begin{ex}\label{ex-computeLCM}
{\rm
For readers unfamiliar with computations using Steinitz numbers we provide an example computation of $LCM(a,b)$. Suppose $a$ and $b$ are Steinitz numbers. Then $a = \prod_{p \in \pi} p^{n(p)}$ and $b = \prod_{p \in \pi} p^{m(p)}$, where $\pi$ is the set of distinct prime numbers. Then 
  $$LCM(a,b) = \prod_{p \in \pi} p^{\max\{n(p),m(p)\}}.$$
In particular, if $\{m_\ell\}_{\ell \geq 1}$ is a sequence of integers, then $LCM\{m_1 \cdot m_2 \cdots m_\ell \mid 1 \leq \ell \leq k\} = m_1 \cdots m_k$, considered as a Steinitz number.  Then $LCM\{m_1 \cdots m_\ell \mid \ell \geq 1\} = \prod_{p \in \pi} p^{n(p)}$ is a Steinitz number, where for each $p \in \pi$ the exponent $n(p)$ is the number of times which $p$ appears as a divisor of the elements in $\{m_\ell \mid \ell \geq 1\}$.
}
\end{ex}

We also note   the   profinite version of Lagrange's Theorem:
\begin{prop}\cite[Proposition~2.1.2]{Wilson1998} \label{prop-lagrange}
 Let $\fK \subset \fH \subset \fG$ be a closed subgroups  of the profinite group $\fG$. Then  
  \begin{equation}\label{eq-relativeorder2}
\Pi[\fG : \fK] =   \Pi[\fG : \fH] \cdot \Pi[\fH : \fK] \ ,
\end{equation}
where the multiplication is taken in the sense of Steinitz numbers. 
\end{prop}

 Now let $(\fX, \G, \Phi)$ be a minimal equicontinuous Cantor action, with   basepoint $x \in \fX$. Recall the \emph{Steinitz orders} of the action, as   in  Definition~\ref{def-steinitzorderaction}: 
  \begin{itemize}
\item $\Pi[\fG(\Phi)] =  LCM \{\# \ \fG(\Phi)/\fN  \mid \fN \subset \fG(\Phi)~ \text{open normal subgroup}\}$, 
\item $\Pi[\fD(\Phi)] =  LCM \{\# \ \fD(\Phi,x)/(\fN \cap \fD(\Phi,x)) \mid \fN \subset \fG(\Phi)~ \text{open normal subgroup}\}$, 
\item $\Pi[\fG(\Phi) : \fD(\Phi)] =  LCM \{\# \ \fG(\Phi)/(\fN \cdot \fD(\Phi,x))  \mid  \fN \subset \fG(\Phi)~ \text{open normal subgroup}\} $.
\end{itemize}
We   consider the dependence of these Steinitz orders on the choices made  and the conjugacy class of the action.
First note that the profinite group $\fG(\Phi)$  does not depend on a choice of basepoint, so this also holds for $\Pi[\fG(\Phi)]$.
 
 Given basepoints $x,y \in \fX$ there exists $\whg_{x,y} \in \fG(\Phi)$ such that $\whg_{x,y}   x = y$. Then the conjugation action of $\whg_{x,y}$ on $\fG(\Phi)$ induces a topological isomorphism of $\fD(\Phi,x)$ with $\fD(\Phi, y)$, and maps a clopen subset of $\fG(\Phi)$ to a clopen subset of $\fG(\Phi)$. Then from the definition, we have    $\Pi[\fD(\Phi,x)] = \Pi[\fD(\Phi, y)]$, and  $\Pi[\fG(\Phi) : \fD(\Phi, x)]  =  \Pi[\fG(\Phi) : \fD(\Phi, y)]$.
 
Let  $(\fX_1, \G_1, \Phi_1)$ and $(\fX_2, \G_2, \Phi_2)$ be isomorphic minimal equicontinuous Cantor actions.  By 
Definition~\ref{def-isomorphism}
there is a homeomorphism $h \colon \fX_1 \to \fX_2$ and group isomorphism $\Theta \colon \G_1 \to \G_2$ so that 
\begin{equation}\label{eq-isomorphism33}
\Phi_1(g) = h^{-1} \circ \Phi_2(\Theta(g)) \circ h   \in   \Homeo(\fX_1) \   \textrm{for  all} \ g \in \G_1 \ .
\end{equation}
Let $\Phi_2' = \Phi_2 \circ \Theta \colon \G_1 \to \Homeo(\fX_2)$, then the images are equal, $\Phi_2(\G) = \Phi_2'(\G)$ and hence so also their closures, $\fG(\Phi_2) = \fG(\Phi_2')$.
The identity  \eqref{eq-isomorphism33} implies that $h$ induces a topological isomorphism between $\fG(\Phi_1)$ and $\fG(\Phi_2')$ and so also between $\fG(\Phi_1)$ and $\fG(\Phi_2)$. 
Thus $\Pi(\fG(\Phi_1)) = \Pi(\fG(\Phi_2))$.

Given $x \in \fX_1$ let $y = h(x) \in \fX_2$, by \eqref{eq-isomorphism33}  the map $h$ induces an isomorphism between $\fD(\Phi_1, x)$ and $\fD(\Phi_2, y)$, and maps clopen subsets of $\fG(\Phi_1)$ to clopen subsets of $\fG(\Phi_2)$. 
Thus  $\Pi[\fD(\Phi_1,x)] = \Pi[\fD(\Phi_2, y)]$ 
and $\Pi[\fG(\Phi_1) : \fD(\Phi_1, x)] =  \Pi[\fG(\Phi_2) : \fD(\Phi_2, y)]$.

 These observations  complete the proof of Theorem~\ref{thm-isoinvariance}.

\subsection{Orders and return equivalence}\label{subsec-re} 
We next consider  how the Steinitz orders behave under return equivalence of actions, and obtain the  proofs of Theorems~\ref{thm-returnequivorder} and \ref{thm-returnequivspectra}.

 Let  $(\fX_1, \G_1, \Phi_1)$ and $(\fX_2, \G_2, \Phi_2)$ be  minimal equicontinuous Cantor actions, and assume that the actions are return equivalent. That is, we assume  there exists 
  an adapted set $U_1 \subset \fX_1$ for the action $\Phi_1$   and  
  an adapted set $U_2 \subset \fX_2$ for the action $\Phi_2$,
such that   the  restricted actions $(U_1, H_{1,U_1}, \Phi_{1,U_1})$ and $(U_2, H_{2,U_2}, \Phi_{2,U_2})$ are isomorphic, with the isomorphism induced by a homeomorphism $h \colon U_1 \to U_2$. 
Thus, the profinite closures 
$$\fH_1 = \overline{H_{1,U_1}} \subset \Homeo(U_1) \, \textrm{ and } \, \fH_2 =    \overline{H_{2,U_2}} \subset \Homeo(U_2)$$ 
are isomorphic. Fix a basepoint $x_1 \in \fX_1$ and set $x_2 = h(x_1) \in U_2$, then the  map $h$ induces an isomorphism between the  isotropy subgroups of the restricted actions, $\fD(\Phi_1|U_1 , x_1)$ and $\fD(\Phi_2|U_2 , x_2)$. 

Our first result is that the asymptotic relative Steinitz order is an invariant of return equivalence.
 \begin{prop}\label{prop-relativeinv}
  Let  $(\fX_1, \G_1, \Phi_1)$ and $(\fX_2, \G_2, \Phi_2)$ be  minimal equicontinuous Cantor actions which are return equivalent. Then 
  \begin{equation} \label{eq-asymporders3}
  \Pi_a[\fG(\Phi_1):\fD(\Phi_1)]   =    \Pi_a[\fG(\Phi_2) : \fD(\Phi_2)]   \ .  
\end{equation}
 \end{prop}
 \proof
For $i=1,2$,    consider the   isotropy subgroup of $U_i$   
\begin{equation}
\fG(\Phi_i)_{U_i} = \left\{ \whg \in \fG(\Phi_i) \mid \whPhi_i(\whg)(U_i) = U_i \right\} \ .
\end{equation}
Then $\fG(\Phi_i)_{U_i}$ is a clopen subgroup   in $\fG(\Phi_i)$, so has finite index $m_i = [\fG(\Phi_i) : \fG(\Phi_i)_{U_i}] = [\G_i : \G_{i,U_i}]$.  Note that since for any $\whg \in \fD(\Phi_i,x_i)$ we have $\whg x = x$, it follows that the action of $\whg$ preserves $U_i$, and so $\fD(\Phi_i,x_i) \subset \fG(\Phi_i)_{U_i}$.

The induced map $\whPhi_i|U_i \colon \fG(\Phi_i)_{U_i} \to \fH_i$ is onto, and 
 the kernel  $\fK_i = \ker \{ \whPhi_i|U_i \colon \fG(\Phi_i)_{U_i} \to \fH_i\}$ is a closed subgroup of $\fG(\Phi_i)_{U_i}$ with $\fK_i \subset \fD(\Phi_i , x_i)$,  since every element of $\fK_i$ fixes $x_i$.
 
Let $\fM_i \subset \fH_i$ be an open   subgroup with $\fD(\Phi_i|U_i , x_i) \subset \fM_i$, then $\fN_i = (\whPhi_i|U_i)^{-1}(\fM_i)$ is an open   subgroup of $\fG(\Phi_i)_{U_i}$ with $\fK_i \subset \fD(\Phi_i, x_i) \subset \fN_i$.  Here $\fD(\Phi_i,x_i)$ is the isotropy group of the action of $\fG(\Phi_i)$ on $\fX_i$, and $\fD(\Phi|U_i,x_i)$ is the isotropy subgroup of the action of $\fH_i \subset Homeo(U_i)$ on $U_i$.

Conversely, let $\fN_i \subset \fG(\Phi_i)_{U_i}$ be an open subgroup with  $\fD(\Phi_i,x_i) \subset \fN_i$. 
Then by \cite[Lemma 2.1.2]{RZ2000}, $\fN_i$ is closed with  finite index in  $\fG(\Phi_i)_{U_i}$   and hence also in $\fG(\Phi_i)$, so  it is clopen hence compact. 
Thus the image $\fM_i = \whPhi_i|U_i(\fN_i) \subset \fH_i$ is a closed subgroup of finite index. Then   \cite[Lemma 2.1.2]{RZ2000} implies it is clopen in $\fH_i$, and 
$\fD(\Phi_i|U_i , x_i) \subset \fM_i$. It follows from Definition~\ref{def-steinitzorderaction} that, for $i=1,2$, 
\begin{equation}\label{eq-relindex}
\Pi[\fG(\Phi_i)_{U_i} :\fD(\Phi_i, x_i)] = \Pi[\fH_i : \fD(\Phi_i |U_i , x_i)] \ .
\end{equation}
The homeomorphism $h \colon U_1 \to U_2$ conjugates the actions $(U_1, \fH_1, \whPhi_1)$ and  $(U_2, \fH_2, \whPhi_2)$ so by the results in Section~\ref{subsec-sodefs} we have   for the restricted actions 
  $$\Pi[\fH_1 : \fD(\Phi_1 |U_1 , x_1)] = \Pi[\fH_2 : \fD(\Phi_2 |U_2 , x_2)].$$
 The equality of the asymptotic Steinitz orders in \eqref{eq-asymporders3} then follows.
 \endproof

 Theorem \ref{thm-returnequivorder} follows immediately from Proposition~\ref{prop-relativeinv}. 
 
 The equality \eqref{eq-relindex} is the key to the proof of Proposition~\ref{prop-relativeinv}. This identity is based on the property that the homomorphism  from  $\fG(\Phi_i)_{U_i}$ to $\fH_i$ has kernel  $\fK_i \subset \fD(\Phi_i, x_i)$,  so the contributions to the Steinitz orders $\fG(\Phi_i)_{U_i}$ and $\fD(\Phi_i, x_i)$ from the subgroup $\fK_i$ cancels out in the relative order $\Pi[\fG(\Phi_i)_{U_i} : \fD(\Phi_i, x_i)]$. However,  the absolute Steinitz orders $\Pi[\fG(\Phi_i)_{U_i}]$ and $\Pi[\fD(\Phi_i, x_i)]$ may indeed include a factor coming from the Steinitz order $\Pi[\fK_i]$. Example~5.3 in \cite{HL2020a} illustrates this.

  For actions with trivial discriminant, Proposition~\ref{prop-relativeinv} has the following consequence:

\begin{cor}\label{cor-2}
Let  $(\fX, \G, \Phi)$   be  a minimal equicontinuous Cantor action with trivial discriminant invariant. Then the asymptotic Steinitz order $\Pi_a[\fG(\Phi)]$ is a return equivalence invariant.
\end{cor}

\proof
In the notation of Proposition \ref{prop-relativeinv}, by assumption we have $\fD(\Phi_1, x_1)$ is the trivial group.  For an adapted clopen set $U_1 \subset \fX_1$ with  $x_1 \in U_1$, we have $\fD(\Phi_1 | U_1, x_1)$ is a quotient of $\fD(\Phi_1, x_1)$ hence is  also trivial. Thus, 
\begin{equation}\label{eq-localconjugacy}
\Pi_a[\fG(\Phi_1)] = \Pi_a[\fG(\Phi_1 | U_1)] = \Pi_a[\fG(\Phi_1 | U_1) : \fD(\Phi_1 | U_1, x_1)] \ .
\end{equation}
Let $(\fX_2,\G_2,\Phi_2)$ be return equivalent to $(\fX_1,\G_1,\Phi_1)$, then the restricted actions $(U_1,H_{1,U_1},\Phi_{1,U_1})$ and $(U_2,H_{2,U_2},\Phi_{2,U_2})$ are isomorphic, which induces a topological isomorphism of the discriminant groups $\fD(\Phi_1|U_1 , x_1)$ and $\fD(\Phi_2|U_2 , x_2)$, and implies that $\fD(\Phi_2|U_2,x_2)$ is trivial.  Using this remark,  a formula analogous to \eqref{eq-localconjugacy} for the action $(\fX_2,\G_2,\Phi_2)$, and Proposition~\ref{prop-relativeinv}, we obtain the claim. 
\endproof

Now consider the behavior of the Steinitz orders $\Pi[\fG(\Phi)]$ and $\Pi[\fD(\Phi, x)]$ under return equivalence of actions. The idea is to use the observation that the action of $\fG(\Phi)$ on $\fX$ is effective (by definition) to construct an effective action map of $\fD(\Phi, x)$ which can be related to a similar construction for a return equivalent action, and so obtain a comparison of their Steinitz orders. 
This yields  the proof of Theorem \ref{thm-returnequivspectra}. 

 Let  $(\fX_1, \G_1, \Phi_1)$ and $(\fX_2, \G_2, \Phi_2)$ be  minimal equicontinuous Cantor actions, and assume that the actions are return equivalent: for   
  an adapted set $U_1 \subset \fX_1$ for the action $\Phi_1$   and  
  an adapted set $U_2 \subset \fX_2$ for the action $\Phi_2$,
  there is a homeomorphism  $h \colon U_1 \to U_2$ which conjugates 
  the  restricted actions $(U_1, H_{1,U_1}, \Phi_{1,U_1})$ and $(U_2, H_{2,U_2}, \Phi_{2,U_2})$.
  
For $i=1,2$,   the action of $\fG(\Phi_i)$ on $\fX_i$ is effective, as $\fG(\Phi_i) \subset \Homeo(\fX_i)$.
Recall that 
$$\fH_i  = \overline{H_{i,U_i}} = \left\{ \whPhi_i(\whg) \mid \whg \in  \fG(\Phi_i)_{U_i} \right\}  \subset  \Homeo(U_i)  \ .$$

   Choose representatives $\{ h_{i,j} \in \G_i \mid   1 \leq j \leq m_i\}$ of the cosets of $\G_i/\G_{i,U_i}$ with $h_{i,1}$ the identity element,  and set 
$U_{i,j} = \Phi_i(h_{i,j})(U_i)$. Thus $U_{i,1} = U_i$, and we have a partition $\fX_i = U_{i,1} \cup \cdots \cup U_{i,m_i}$.

Introduce the normal core of $\fG(\Phi_i)_{U_i}$ given by 
\begin{equation}
\fN(\Phi_i) = \bigcap_{j=1}^{m_i} \ \Phi_i(h_{i,j})^{-1} \cdot \fG(\Phi_i)_{U_i} \cdot \Phi_i(h_{i,j}) \subset \fG(\Phi_i)_{U_i} \ ,
\end{equation}
which is a clopen subgroup of $\fG(\Phi_i)$ of finite index $n_i = [\fG(\Phi_i) : \fN(\Phi_i)]$, where $m_i$ divides $n_i$.
In particular, we have $[\fG(\Phi_i)_{U_i} : \fN(\Phi_i)] < n_i$.

The fact that $\fN(\Phi_i)$ is a normal subgroup of $\fG(\Phi_i)$ implies  that the action of $\fN(\Phi_i)$ on the partition of $\fX_i$ maps each of the sets $U_{i,j}$ to itself.
   
 Recall that $\whPhi_i: \fG(\Phi_i) \to Homeo(\fX_i)$ is the action of the profinite completion of $(\fX_i,\G_i,\Phi_i)$, $i=1,2$.   
    
For $\whg \ne \whe$,  the action of $\whPhi_i(\whg)$ on $\fX_i$ is non-trivial, so if  $\whg \in \fN(\Phi_i)$ also, then for some $1 \leq j \leq m_i$ the restricted action of $\whPhi_i(\whg)$ on $U_{i,j}$ must be non-trivial. That is, for some $j$ we have 
\begin{equation}
\whg  \not\in   \ker \left\{ \whPhi_{i,j} \equiv  \whPhi_i | U_{i,j} \colon \fN(\Phi_i) \to \Homeo(U_{i,j})   \right\} \ . 
\end{equation}
Define a representation $\whrho_i$ of $\fN(\Phi_i)$ into a product of $m_i$ copies of $\fH_i$ by setting, for $\whg \in \fN(\Phi_i)$, 
\begin{equation}\label{eq-prodrep}
\whrho_i \colon \fN(\Phi_i) \to \fH_i  \times \cdots \times \fH_i  \quad , \quad \whrho_i(\whg) = \whPhi_i^1(\whg)  \times \cdots \times \whPhi_i^{m_i}(\whg) \ ,
\end{equation}
where we use that $\fN(\Phi_i)$ is normal in $ \fG(\Phi_i)$, so for $\whg \in \fN(\Phi_i)$ the following is well-defined:
$$\whPhi_i^j(\whg) = \Phi_i(h_{i,j})^{-1} \circ \whPhi_{i,j}(\whg)  \circ \Phi_i(h_{i,j}) =   \whPhi_i(h_{i,j}^{-1} \ \whg \ h_{i,j}) | U_{i}  \in \fH_i \ .$$
The kernel of  $\whrho_i$ is trivial by the above arguments,  so there is an isomorphism $\fN(\Phi_i) \cong \whrho_i(\fN(\Phi_i))$.
  This diagonal   trick to obtain the injective map $\whrho_i$ was first used   in the proof of   \cite[Theorem~1.2]{HL2020a}.

  The index $n_i = [\fG(\Phi_i) : \fN(\Phi_i)] < \infty$, so we have     \begin{equation}\label{eq-compare1}
   [\fG(\Phi_i)] \mor [\fG(\Phi_i)_{U_i}] \mor [\fN(\Phi_i)] = [\whrho_i(\fN(\Phi_i))] \ .
   \end{equation}

 Let $p_{i,1} \colon \fH_i \times \cdots \times \fH_i \to \fH_i$ denote the projection onto the first factor.
    Then the composition   $p_{i,1} \circ \whrho_i$  equals  the restriction to $\fN(\Phi_i)$ of the map
  $\whPhi_{i,U_i} \colon \fG(\Phi_i)_{U_i} \to \fH_i$. Let 
  $$\fL_i = \ker \ p_{i,1} \colon  \whrho_i(\fN(\Phi_i)) \to \whPhi_{i,U_i}(\fN(\Phi_i))$$ 
  denote the kernel of the restriction of $p_{i,1}$. Then by 
    Proposition~\ref{prop-lagrange} applied to the inclusions $\{\whe\} \subset \fL_i \subset \whrho_i(\fN(\Phi_i))$, by the 
identity \eqref{eq-relativeorder2} we have 
$  \Pi[\whrho_i(\fN(\Phi_i))] =   \Pi[\whrho_i(\fN(\Phi_i)) : \fL_i] \cdot \Pi[\fL_i ]$.

Since by the first isomorphism theorem $\whPhi_{i,U_i}(\fN(\Phi_i)) = \whrho_i(\fN(\Phi_i)) /\fL_i$, then
$$ \Pi[\whrho_i(\fN(\Phi_i)) : \fL_i] = \Pi[\whPhi_{i,U_i}(\fN(\Phi_i))],$$
and thus we have the inequality of Steinitz orders $ [\whPhi_{i,U_i}(\fN(\Phi_i))]  \leq [\whrho_i(\fN(\Phi_i))]$.

Now note that     $\fN(\Phi_i)$  has finite index in $\fG(\Phi_i)_{U_i}$ implies the same holds for its image under $\whPhi_{i,U_i}$, so we have
  $[\whPhi_{i,U_i}(\fN(\Phi_i))] \mor [\fH_i]$. 
  Thus we have the   estimate on Steinitz orders
  \begin{equation}\label{eq-compare2}
[\fH_i] \mor [\whPhi_{i,U_i}(\fN(\Phi_i))]  \leq [\whrho_i(\fN(\Phi_i))] \ .
\end{equation}

On the other hand, from the embedding in \eqref{eq-prodrep} we have 
\begin{equation}\label{eq-compare3}
[\whrho_i(\fN(\Phi_i))] \leq [\fH_i] \cdots [\fH_i] = [\fH_i]^{m_i} \ . 
\end{equation}
Combining the estimates \eqref{eq-compare1}, \eqref{eq-compare2}, and \eqref{eq-compare3} we obtain that $\pi_{\infty}([\fH_i]) = \pi_{\infty}(\fN(\Phi_i)) = \pi_{\infty}(\fG(\Phi_i))$. Moreover,  $\pi_{f}([\fH_i])$  and  $\pi_{f}(\fG(\Phi_i))$ differ by at most a finite subset of primes. 
As $\fH_1$ and $\fH_2$ are topologically isomorphic, this shows that the prime spectra of $\fG(\Phi_1)$  and $\fG(\Phi_2)$
satisfy the claim of Theorem~\ref{thm-returnequivspectra}.  

We can apply the same analysis as above to the isotropy subgroups $\fD(\Phi_1, x_1)$  and $\fD(\Phi_2, x_2)$ to obtain the stated relations between their prime spectra, completing the proof of  Theorem~\ref{thm-returnequivspectra}.

\subsection{Algebraic model} \label{subsec-algmodel} 

In this section we reformulate the abstract  Definition~\ref{def-steinitzprofinite} of the Steinitz order invariants in terms of  an  algebraic model for a Cantor action. This  provides an effective method of calculating and working with these invariants.
  We first  recall the construction of the algebraic models for an action  $(\fX, \G, \Phi)$ and  its profinite completion.
 
 For $x \in \fX$, by Proposition~\ref{prop-adpatedchain} there exists an adapted neighborhood basis $\cU = \{U_{\ell} \subset \fX  \mid \ell \geq 0\}$ at $x$ for the action $\Phi$.
  Let    $\G_{\ell} = \G_{U_{\ell}}$ denote the stabilizer group of $U_{\ell}$.
Then  we obtain a strictly descending chain of finite index subgroups 
 \begin{equation}\label{eq-groupchain}
 \cG^x_{\cU} = \{\G = \G_0 \supset \G_1 \supset \G_2 \supset \cdots \} \ .
\end{equation}
Note that each $\G_{\ell}$ has finite index in $\G$, and is not assumed to be a normal subgroup.  Also note that while the intersection of the chain $\cU$ is a single point $\{x\}$, the intersection of the stabilizer groups   in  $\cG^x_{\cU}$ need not be the trivial group.
 
Next, set $X_{\ell} = \G/\G_{\ell}$ and note that  $\G$ acts transitively on the left on   $X_{\ell}$.    
The inclusion $\G_{\ell +1} \subset \G_{\ell}$ induces a natural $\G$-invariant quotient map $p_{\ell +1} \colon X_{\ell +1} \to X_{\ell}$.
 Introduce the inverse limit 
 \begin{eqnarray} 
X_{\infty} & \equiv &  \lim_{\longleftarrow} ~ \{ p_{\ell +1} \colon X_{\ell +1} \to X_{\ell}  \mid \ell \geq 0 \} \label{eq-invlimspace}\\
& = &  \{(x_0, x_1, \ldots ) \in X_{\infty}  \mid p_{\ell +1 }(x_{\ell + 1}) =  x_{\ell} ~ {\rm for ~ all} ~ \ell \geq 0 ~\} ~ \subset \prod_{\ell \geq 0} ~ X_{\ell} \  ,  \nonumber
\end{eqnarray}
which is a Cantor space with the Tychonoff topology, and the actions of $\G$ on the factors $X_{\ell}$ induce    a minimal  equicontinuous action on the inverse limit, denoted by  $\Phi_x \colon G \times X_{\infty} \to X_{\infty}$. Denote the points in $X_{\infty}$ by 
$x = (x_{\ell}) \in X_{\infty}$. There is a natural basepoint $x_{\infty} \in X_{\infty}$ given by the cosets of the identity element $e \in \G$, so $x_{\infty} = (e \G_{\ell})$. A basis of  neighborhoods of $x_{\infty}$ is given by the clopen sets 
\begin{equation}\label{eq-openbasis}
U_{\ell} = \left\{ x = (x_{\ell}) \in X_{\infty}   \mid  x_i = e \G_i \in X_i~, ~ 0 \leq i < \ell ~  \right\} \subset X_{\infty} \ .
\end{equation}

 For each $\ell \geq 0$, we have the ``partition coding map'' $\Theta_{\ell} \colon \fX \to X_{\ell}$ which is $G$-equivariant.  The maps $\{\Theta_{\ell}\}$ are compatible with the   map on quotients in \eqref{eq-invlimspace}, and so they induce a  limit map $\Theta_x \colon \fX \to X_{\infty}$. The fact that the diameters of the clopen sets $\{U_{\ell}\}$ tend to zero, implies that $\Theta_x$ is a homeomorphism.  Moreover, $\Theta_x(x) =  x_{\infty} \in X_{\infty}$.  
 \begin{thm}\cite[Appendix~A]{DHL2016a}
 The map $\Theta_x \colon  \fX \to X_{\infty}$ induces an isomorphism of the Cantor actions $(\fX,\G,\Phi)$ and $(X_{\infty}, \G, \Phi_x)$.
 \end{thm}
 The   action $(X_{\infty}, G, \Phi_x)$   is called the \emph{odometer model} centered at $x$ for the action $(\fX,\G,\Phi)$.
 The dependence of the model   on the choices of a base point $x \in \fX$ and adapted neighborhood basis $\cU$ is discussed in detail in the   works \cite{DHL2016a,FO2002,HL2018a,HL2019a}.

 Next, we develop the algebraic model for the profinite action  $\whPhi \colon \fG(\Phi) \times \fX \to \fX$ of the completion  $\fG(\Phi) \equiv \overline{\Phi(\G)} \subset \Homeo(\fX)$. 
Fix a choice of group chain $\{\G_{\ell} \mid \ell \geq 0\}$ as above, which provides an algebraic model for the action $(\fX,\G,\Phi)$.

 For each $\ell \geq 1$, let $C_{\ell} \subset \G_{\ell}$ denote the \emph{core} of   $\G_{\ell}$, that is, the largest normal subgroup of   $\G_{\ell}$. So 
\begin{equation}\label{eq-core}
C_{\ell} ~ = {\rm Core}(\G_{\ell}) ~ = ~ \bigcap_{g \in \G} ~ g \ \G_{\ell} \ g^{-1} ~ \subset \G_{\ell} ~ .
\end{equation}
As   $\G_{\ell}$ has finite index in $\G$, the same holds for $C_{\ell}$. Observe that for all $\ell \geq 0$,   we have $C_{\ell +1} \subset C_{\ell}$.

Introduce the quotient group  $Q_{\ell} = \G/C_{\ell}$ with identity element $e_{\ell} \in Q_{\ell}$. There are natural quotient maps $q_{\ell+1} \colon Q_{\ell +1} \to Q_{\ell}$, and we can form the inverse limit group
  \begin{eqnarray} 
\whGamma_{\infty} & \equiv &  \lim_{\longleftarrow} ~ \{ q_{\ell +1} \colon Q_{\ell +1} \to Q_{\ell}  \mid \ell \geq 0 \} \label{eq-invgroup}\\
& = &  \{(g_{\ell}) = (g_0, g_1, \ldots )    \mid g_{\ell} \in Q_{\ell} ~ , ~ q_{\ell +1 }(g_{\ell + 1}) =  g_{\ell} ~ {\rm for ~ all} ~ \ell \geq 0 ~\} ~ \subset \prod_{\ell \geq 0} ~ \G_{\ell} ~ , \label{eq-coordinates}
\end{eqnarray}
which is a Cantor space with the Tychonoff topology. The left actions of $\G$ on the spaces $X_{\ell} = \G/\G_{\ell}$ induce    a minimal  equicontinuous action of $\whGamma_{\infty}$ on $X_{\infty}$, again denoted by  $\whPhi \colon \whGamma_{\infty} \times X_{\infty} \to X_{\infty}$. Note that the isotropy group of the identity coset of the action of $Q_{\ell} = \G_{\ell}/C_{\ell}$ on $X_{\ell}=  \G/\G_{\ell}$ is the subgroup $D_{\ell} = \Gamma_{\ell}/C_{\ell}$.

Denote the points in $\whGamma_{\infty}$ by 
$\whg = (g_{\ell}) \in \whGamma_{\infty}$ where $g_{\ell} \in Q_{\ell}$. There is a natural basepoint $\whe_{\infty} \in \whGamma_{\infty}$ given by the cosets of the identity element $e \in \G$, so $\whe_{\infty} = (e_{\ell})$ where $e_{\ell} = e C_{\ell} \in Q_{\ell}$ is the identity element in $Q_{\ell}$. 

For each $\ell \geq 0$, let $\Pi_{\ell} \colon \whGamma_{\infty} \to Q_{\ell}$ denote the projection onto the $\ell$-th factor in \eqref{eq-invgroup}, so in the coordinates of \eqref{eq-coordinates}, we have $\Pi_{\ell}(\whg) = g_{\ell} \in Q_{\ell}$. 
The maps $\Pi_{\ell}$ are continuous for the profinite topology on $\whGamma_{\infty}$, so the pre-images of points in $Q_{\ell}$ are clopen subsets. In particular, the  fiber of $Q_{\ell}$  over   $e_{\ell}$ is the normal subgroup
\begin{equation}\label{eq-opennbhds}
\whC_{\ell} = \Pi_{\ell}^{-1}(e_{\ell}) =  \{(g_{\ell})  \in \whGamma_{\infty}  \mid  g_{i} \in C_{i} ~ , ~ 0 \leq i \leq \ell \} \ . 
\end{equation}

Then the collection $\{\whC_{\ell} \mid \ell \geq 1\}$ forms a basis of   clopen neighborhoods of $\whe_{\infty} \in \whGamma_{\infty}$. That is, for each clopen set $\whU \subset \whGamma_{\infty}$ with $\whe_{\infty} \in \whU$, there exists $\ell_0 > 0$ such that $\whC_{\ell} \subset \whU$ for all $\ell \geq \ell_0$.
\begin{thm}\cite[Theorem~4.4]{DHL2016a}\label{thm-fundamentaliso}
There is an   isomorphism $\whtau \colon \fG(\Phi) \to \whGamma_{\infty}$ which conjugates the profinite action
$(\fX, \fG(\Phi), \whPhi)$ with the profinite action
$(X_{\infty}, \whGamma_{\infty}, \whPhi)$. In particular, $\whtau$ 
identifies the isotropy group $\fD(\Phi, x) = \fG(\Phi)_{x}$ with the inverse limit subgroup
\begin{equation}\label{eq-discformula}
D_{\infty}  = \varprojlim \ \{q_{\ell +1} \colon \G_{\ell +1}/C_{\ell+1} \to \G_{\ell}/C_{\ell} \mid \ell \geq 0\} \subset \whGamma_{\infty}~ .
\end{equation}
\end{thm}

The maps $q_{\ell +1}$ in the formula \eqref{eq-discformula} need not be surjections, and thus the calculation of the inverse limit $D_{\infty}$   can involve some subtleties. For example, it is possible that each group $Q_{\ell}$ is non-trivial for $\ell > 0$, and yet $D_{\infty}$ is the trivial group (see Example~\ref{ex-trivial}.) This phenomenon leads to the following considerations.
Observe that the formula \eqref{eq-discformula} implies  the restriction of the projection map  $\Pi_{\ell} \colon D_{\infty} \to Q_\ell$ yields a map $\Pi_{\ell} \colon D_{\infty} \to D_{\ell} \equiv \G_{\ell}/C_{\ell} \subset Q_{\ell}$.
Set
\begin{equation}\label{eq-discimage}
D_{\ell}^* =  \Pi_{\ell}(D_{\infty})  \subset D_{\ell} \ .
\end{equation}
 We recall a concept definition from \cite[Definition~5.6]{DHL2016a}:
 
 \begin{defn}\label{def-normalform}
A group chain    $\{\G_{\ell} \mid \ell   \geq 0\}$   in $\G$ is in \emph{normal form} if    $D_{\ell}^* = D_{\ell}$, for   $\ell \geq 0$.
\end{defn}
 Recall that if the group chain $\{\G_{\ell} \mid \ell   \geq 0\}$ is in normal form, 
 then each of the bonding maps $q_{\ell +1}$ in \eqref{eq-discformula} is a surjection.  We note that, given any group chain $\cG = \{\G_\ell \mid \ell \geq 0\}$, by \cite[Proposition~5.7]{DHL2016a} there exists a group chain  $\cG' = \{\G_{\ell}' \mid \ell \geq 0\}$    in normal form which is equivalent to $\cG$, that is, up to a choice of infinite subsequences the group chains are intertwined, $\G_0 \supset \G_1 ' \supset \G_1 \supset \G_2' \supset \cdots$. As explained in \cite{DHL2016a}, the actions defined by equivalent group chains $\cG$ and $\cG'$ using formulas \eqref{eq-invlimspace} - \eqref{eq-invgroup} are  isomorphic, and the   homeomorphism implementing the isomorphism preserves the basepoint.

 \subsection{Steinitz orders for algebraic models} \label{subsec-steinitzalgmodel}

 Let $(\fX, \G, \Phi)$ be   a minimal equicontinuous Cantor action, chose $x \in \fX$ and an adapted neighborhood basis $\cU$ at $x$, then let $\cG = \{\G_{\ell} \mid \ell \geq 0\}$ be the associated group chain formed by the stabilizer subgroups of the clopen sets $U_{\ell}$ in $\cU$. We continue  further with the notation in Section~\ref{subsec-algmodel}.

For $\ell \geq 0$, we have   the finite sets $X_{\ell} = \G/\G_{\ell}$, and the finite groups  $Q_{\ell} = \G/C_{\ell}$, 
$D_{\ell} = \G_{\ell}/C_{\ell}$ and $D_{\ell}^* =  \Pi_{\ell}(\cD_{\infty}) \subset D_{\ell}$.
Introduce the sequences of integers:
\begin{equation}\label{eq-dims}
m_{\ell} = \# \ X_{\ell}  \quad ; \quad n_{\ell} = \# \ Q_{\ell} \quad ; \quad k_{\ell} = \# \ D_{\ell} \quad ; \quad k_{\ell}^* = \# \ D_{\ell}^* \ .
\end{equation} 
We make some elementary observations about these sequences of integers. 

Lagrange's Theorem implies that $n_{\ell} = m_{\ell} k_{\ell}$  for $\ell \geq 0$, and we also have   $k_{\ell}^* \leq k_{\ell}$. 

Note that $m_{\ell +1} = m_{\ell} \cdot [\G_{\ell} : \G_{\ell+1}]$. As   the inclusion $\G_{\ell +1} \subset \G_{\ell}$ is proper, we have $[\G_{\ell} : \G_{\ell+1}] > 1$ and so  $\{ m_{\ell} \mid \ell \geq 0\}$ is a strictly increasing sequence.

Also, $C_{\ell +1} \subset C_{\ell}$,  and $n_{\ell+1} = n_{\ell} \cdot [C_{\ell} : C_{\ell+1}]$ so  $\{ n_{\ell} \mid \ell \geq 0\}$ is a non-decreasing sequence.

As $k_{\ell}^*$ is the order of the projection of  $\cD_{\infty}$ into $Q_{\ell}$, the sequence $\{ k_{\ell}^* \mid \ell \geq 0\}$ is   non-decreasing.  For instance, when $\cD_\infty$ is a finite group, then there exist $m \geq 0$ such that $k_\ell^* = k_{\ell+1}^*$ for all $\ell \geq m$.

 \begin{prop}\label{prop-steinitzinvariance}
  Let $(\fX, \G, \Phi)$ be   a minimal equicontinuous Cantor action.
Given a basepoint  $x \in \fX$, and an adapted neighborhood basis $\cU$ at $x$,   let $\cG = \{\G_{\ell} \mid \ell \geq 0\}$ be the associated group chain formed by the stabilizer subgroups of the clopen sets $U_{\ell}$ in $\cU$. Then the Steinitz orders for the action, as defined in Definition~\ref{def-steinitzorderaction},  can be calculated as follows:
\begin{enumerate}
\item \quad  $\Pi[\fG(\Phi)] = LCM \ \{n_{\ell} \mid \ell \geq 0\}$, 
\item \quad  $\Pi[\fG(\Phi) : \fD(\Phi , x)] = LCM \ \{m_{\ell} \mid \ell \geq 0\}$,  
\item \quad  $\Pi[\fD(\Phi, x)] = LCM \ \{k^*_{\ell}  \mid \ell \geq 0\} \leq LCM \ \{k_{\ell}  \mid \ell \geq 0\}$ \ .
\end{enumerate}
\end{prop}
\proof
By Theorem~\ref{thm-fundamentaliso}, 
there is an isomorphism $\whtau \colon \fG(\Phi) \to \whGamma_{\infty}$ 
which conjugates the profinite action $(\fX, \whGamma, \whPhi)$
  with the profinite action  $(X_{\infty}, \whGamma_{\infty}, \whPhi)$. 
 By the results of Section~\ref{subsec-sodefs}, it suffices to show  that the formulas in Proposition \ref{prop-steinitzinvariance}, (1)-(3),  hold for the action $(X_{\infty}, \whGamma_{\infty}, \whPhi)$.

Recall that $\whC_{\ell}$ is the normal clopen subgroup of $\whGamma_{\infty}$ defined  in  \eqref{eq-opennbhds}. Since  $\{\whC_\ell\}_{\ell \geq 0}$  form a neighborhood basis for the identity in  $\whGamma_\infty$, for any clopen normal subgroup $\cN \subset \whGamma_{\infty}$, there exists $\ell > 0$ such that $\whC_{\ell} \subset \cN$. It follows that
$\# (\whGamma_{\infty}/\cN)$ divides $\# (\whGamma_{\infty}/\whC_{\ell}) = \# Q_{\ell}$. Noting that $\whC_{\ell}$ is itself a clopen normal subgroup, we have
\begin{eqnarray}
\lefteqn{ LCM \{\# \  \whGamma_{\infty}/\cN   \mid \cN \subset \whGamma_{\infty} ~ \text{clopen normal subgroup} \}     = }  \label{eq-reduction}  \\
 & &   LCM \{\# \  \whGamma_{\infty}/\whC_{\ell}     \mid  \ell > 0\} =   LCM \{\# \  Q_{\ell}  \mid  \ell > 0\} \ . \nonumber
\end{eqnarray}
Then by  Definition \ref{def-steinitzorderaction},  
\begin{eqnarray*}
\Pi[\fG(\Phi)] & = &   LCM \{\# \  \fG(\Phi)/\fN   \mid \fN \subset \fG(\Phi)~ \text{clopen normal subgroup}\}  \\
  & = &   LCM \{\# \  \whGamma_{\infty}/\cN   \mid \cN \subset \whGamma_{\infty} ~ \text{clopen normal subgroup} \}  \\
  & = &     LCM \{\# \  Q_{\ell}  \mid  \ell > 0\}   =   LCM \{ \ n_{\ell}  \mid  \ell > 0\}   \ .
\end{eqnarray*}

The proofs of   the identities (2) and (3) in Proposition \ref{prop-steinitzinvariance}  require an additional consideration. Introduce the closures of the subgroups $\G_{\ell}$,  for $\ell > 0$,
\begin{equation}
\whGamma_{\ell} = \overline{\G_{\ell}}  = \left\{ \whg = (g_{\ell}) \in \whGamma_{\infty}  \mid g_i = e_i ~, ~ 0 \leq i < \ell ~ ; ~ g_i \in \G_i ~ , ~ i \geq \ell \right\}   \subset \whGamma_{\infty} \ .
\end{equation}
Then each  $\whGamma_{\ell}$ is a    clopen subset  of $\whGamma_{\infty}$, and from the formula \eqref{eq-discformula} we have $D_{\infty} \subset \whGamma_{\ell}$ for all $\ell \geq 0$, and moreover, we have
\begin{equation}\label{eq-intersection}
D_{\infty} = \bigcap_{\ell > 0} ~ \whGamma_{\ell} \ .
\end{equation}
The equality in \eqref{eq-intersection} follows as the action of $\whg \in \Gamma_{\ell}$ on $X_{\infty}$ fixes  the clopen set $U_{\ell}$ defined by \eqref{eq-openbasis}, so $\whg \in \whGamma_{\ell}$ for all $ \ell > 0$ implies that its action on $X_{\infty}$ fixes the intersection, $x_{\infty} = \cap_{\ell > 0} U_{\ell}$.
Also, observe that for $\ell > 0$ we have the identity
\begin{equation}
\G_{\ell} = \left\{ g \in \G \mid \whg = (g,g, \ldots) \in  \whGamma_{\ell} \right\} \ ,
\end{equation}
and consequently there is an isomorphism $\whGamma_{\infty}/\whGamma_{\ell} \cong \G/\G_{\ell}$.

Next,  observe that given a clopen normal subgroup $\cN \subset \whGamma_{\infty}$, by \eqref{eq-intersection} there exists $\ell$ such that $\whGamma_{\ell} \subset   \cN \cdot D_{\infty}$. For instance, this holds for any $\ell \geq 0$ such that $\whC_\ell \subset \cN$. Then the identity (2) in  Proposition \ref{prop-steinitzinvariance}  follows from the fact that $\whGamma_{\ell}$ is a clopen neighborhood of $D_{\infty}$, and reasoning as for \eqref{eq-reduction}, we have  
  \begin{eqnarray*}
\Pi[\fG(\Phi)  : \fD(\Phi , x)] & = &   LCM \{\# \  \fG(\Phi)/(\fN \cdot \fD(\Phi , x))   \mid \fN \subset \fG(\Phi)~ \text{clopen normal subgroup}\}  \\
& = &   LCM \{\# \  \whGamma_{\infty}/(\cN \cdot D_{\infty})   \mid \cN \subset \whGamma_{\infty}~ \text{clopen normal subgroup}\}  \\
  & = &   LCM \{\# \  \whGamma_{\infty}/\whGamma_{\ell}    \mid  \ell > 0\}   =  LCM \{\# \  \G/\G_{\ell}    \mid  \ell > 0\}    \\
  & = &    LCM \{ \ m_{\ell}  \mid  \ell > 0\}   \ .
\end{eqnarray*}

Similarly, the proof of  the identity (3) in   Proposition \ref{prop-steinitzinvariance}  follows from the calculations:
\begin{eqnarray*}
\Pi[\fD(\Phi , x)] & = &   LCM \{\# \  \fD(\Phi , x)/(\fN \cap \fD(\Phi , x))  \mid \fN \subset \fG(\Phi)~ \text{clopen normal subgroup}\}  \\
  & = &   LCM \{\# \  D_{\infty}/(\cN  \cap  D_{\infty})  \mid \cN \subset \whGamma_{\infty} ~ \text{clopen normal subgroup} \}  \\
  & = &   LCM \{\# \ D_{\infty}/  (\whC_{\ell}    \cap  D_{\infty})   \mid  \ell > 0\}  \\ 
   & = &   LCM \{\# \ \Pi_{\ell}(\cD_{\infty})   \mid  \ell > 0\}  
    =    LCM \{\# \ k^*_{\ell}   \mid  \ell > 0\}  \ . 
\end{eqnarray*}
This completes the proof of Proposition~\ref{prop-steinitzinvariance}.
\endproof

 As remarked is the discussion of Definition~\ref{def-normalform}, the  condition that the chain $\cG$ is descending  does not impose sufficient restrictions on the behavior of the orders of the groups $D_{\ell} = \G_\ell/C_\ell$ in order to compute $\Pi[\fD(\Phi , x)]$.
 Rather,  computing $LCM\{D_\ell \mid \ell \geq 0\} = LCM\{k_\ell \mid \ell \geq 0\}$  yields an upper bound on the Steinitz order of $D_{\infty}$. However, if we are given  that the chain $\cG$ is in normal form, as in Definition~\ref{def-normalform},   then this indeterminacy is removed.
 \begin{cor}
   Let $\cG = \{\G_{\ell} \mid \ell \geq 0\}$ be a group chain in normal form which gives an algebraic model for   a Cantor action  $(\fX, \G, \Phi)$. Then we have
 \begin{equation}\label{eq-SOnormalform}
\Pi[\fD(\Phi , x)] =   LCM \{\# \ D_{\ell}   \mid  \ell > 0\}  =  LCM \{\# \ k_{\ell}   \mid  \ell > 0\} \ . 
\end{equation}
  \end{cor}
 It is often the case when constructing examples of Cantor actions, that the normal form property is guaranteed by the choices in the construction, and then \eqref{eq-SOnormalform} calculates the Steinitz order of the discriminant of the action.

\subsection{Steinitz orders of solenoidal manifolds}\label{subsec-solenoids}
 
 We   relate the asymptotic Steinitz order for a tower of coverings with the Steinitz order invariants for Cantor actions.
This yields the proof of Theorem~\ref{thm-main1}.   We first recall some preliminary constructions for solenoidal manifolds.

Let $M_0$ be a compact connected manifold without boundary. Let 
$\cP = \{\, q_{\ell} \colon  M_{\ell} \to M_{\ell -1} \mid  \ell \geq 1\}$ be a presentation as in Section~\ref{sec-intro}.  Let $\cS_{\cP}$ be the inverse limit    of this presentation as in \eqref{eq-presentationinvlim}. A point $x \in \cS_{\cP}$ is represented by a sequence, $x = (x_0, x_1, \ldots )$ with $x_{\ell} \in M_{\ell}$. 
For each $\ell \geq 0$, projection onto the $\ell$-th factor in \eqref{eq-presentationinvlim} yields a fibration denoted by 
  $\whq_{\ell} \colon \cS_{\cP} \to M_{\ell}$, so $\whq_{\ell}(x) = x_{\ell}$.   
Denote the iterated covering map by
$\ovq_{\ell} = q_{\ell} \circ q_{\ell -1} \circ \cdots \circ q_1 \colon M_{\ell} \to M_0$, and note that $\whq_0 = \ovq_{\ell} \circ \whq_{\ell}$.

Choose a basepoint $x_0 \in M_0$, and let $\fX_0 = \whq_0^{-1}(x_0)$ denote the fiber of the projection map $\whq_0$.
Then $\fX_0$ is a Cantor space,  and  the holonomy along the leaves of the foliation $\F_{\cP}$ on $\cS_{\cP}$ induces 
the \emph{monodromy action} of  the fundamental group $\G_0 = \pi_1(M_0, x_0)$     on $\fX_0$. This action is  discussed in greater detail is many works, for example in  \cite{CandelConlon2000}.

Choose a basepoint $x  \in \fX_0$ and then for each $\ell \geq 0$, set $x_{\ell} = \whq_{\ell}(x) \in M_{\ell}$.
Then $\ovq_{\ell}(x_{\ell}) = x_0$ so we get induced maps of fundamental groups, 
$(\ovq_{\ell})_{\#} \colon \pi_1(M_{\ell} , x_{\ell}) \to \pi_1(M, x_0) = \G_0$. Let $\G_{\ell} \subset \G_0$ denote the image of this map, so $\G_{\ell} \subset \G_0$  is a subgroup of finite index. Note that $\ovq_{\ell} \colon M_{\ell} \to M_0$ is a normal covering map exactly when $\G_{\ell}$ is a normal subgroup of $\G_0$.

Let $(X_{\infty}, \G_0, \Phi_x)$ be the Cantor action associated to the group chain $\cG_{x} = \{\G_{\ell} \mid \ell \geq 0 \}$ constructed in Section~\ref{subsec-algmodel} above. Then the monodromy action of $\G_0$ on $\fX_0$ determined by the foliation on  $\cS_{\cP}$  is conjugate to the action 
$(X_{\infty}, \G, \Phi_x)$, as discussed in \cite[Section~2]{DHL2016b} and \cite[Section~3.1]{DHL2016c}.
In particular, note that the degree  of the covering map $\ovq_{\ell} \colon M_{\ell} \to M_0$ equals the index $\#[\G_0 : \G_{\ell}]$.
Thus,  by the identity (2) in Proposition~\ref{prop-steinitzinvariance},  the {Steinitz order} $\Pi[\cP] $ of   $\cP$ in Definition~\ref{def-steinitzpres} equals the relative Steinitz order $\Pi[\fG(\Phi_x) : \fD(\Phi_x , x)]$ of the action $(X_{\infty}, \G, \Phi_x)$.
 
  Now suppose, for $i=1,2$,  we are given a solenoidal manifold $\cS_{\cP_i}$ defined by the presentation $\cP_i$ and there exists a homeomorphism   $h \colon \cS_{\cP_1} \to \cS_{\cP_2}$. Then by the results of Section~\ref{subsec-morita},   the  homeomorphism $h$ induces a return equivalence of their monodromy actions, and thus the algebraic models for these actions defined by $\cP_1$ and $\cP_2$ are return equivalent. 
  
   By Proposition~\ref{prop-relativeinv} we have 
   $\Pi_a[\fG(\Phi_1):\fD(\Phi_1)]   =    \Pi_a[\fG(\Phi_2) : \fD(\Phi_2)]$.
 
 Proposition~\ref{prop-steinitzinvariance} identifies  $\Pi_a[\fG(\Phi_i):\fD(\Phi_i)]$ with the asymptotic Steinitz order $\Pi_a[\cP_i]$ and so we obtain the conclusion of Theorem~\ref{thm-main1}.

 \section{Nilpotent actions}\label{sec-nilpotent}

In this section, we apply the notion of  the Steinitz order of a nilpotent Cantor action to the study of  its dynamical properties. 
The proof of Theorem~\ref{thm-nilstable}    is based on  the special properties of the profinite completions of nilpotent groups, 
in particular the uniqueness of their Sylow $p$-subgroups, and the relation of this algebraic property with the dynamics of the action.

\subsection{Noetherian groups}\label{subsec-noetherian}
Baer introduced the notion of a Noetherian group in his work \cite{Baer1956}. A countable group $\G$ is said to be \emph{Noetherian} if every increasing chain of   subgroups $\{H_i \mid i \geq 1 \}$ of $\G$ has a maximal element $H_{i_0}$. Equivalently, $\G$ is Noetherian if every increasing chain of subgroups in $\G$ eventually stabilizes.  
It is easy to see that the group $\mZ$ is Noetherian, that a finite product of Noetherian groups is Noetherian, and that a subgroup and quotient group of a Noetherian group is Noetherian.  Thus, a finitely-generated nilpotent group is Noetherian. 
  
The notion of a Noetherian group has a generalization which is useful for the study of actions of  profinite groups (see \cite[page 153]{Wilson1998}.)  
\begin{defn} \label{def-noetherian} 
A profinite group $\fG$ is said to be \emph{topologically Noetherian} if every increasing chain of \emph{closed} subgroups $\{\fH_i \mid i \geq 1 \}$ of $\fG$ has a maximal element $\fH_{i_0}$.
\end{defn}
   
   We illustrate this concept with two canonical examples of profinite completions of $\mZ$. First, let $\whmZ_p$ denote the $p$-adic integers, for $p$ a prime. That is, $\whmZ_p$ is the completion of $\mZ$ with respect to the chain of subgroups 
 $\cG = \{\G_{\ell} = p^{\ell} \mZ \mid \ell \geq 1\}$. The closed subgroups of  $\whmZ_p$ are given by   $p^i \cdot \whmZ_p$ for some fixed $i > 0$, hence satisfy the ascending chain property  in Definition \ref{def-noetherian}. 
 
 Next, let $\pi = \{p_i \mid i \geq 1\}$ be an infinite collection of distinct primes, and define the increasing chain of subgroups of $\mZ$ defined by  $\cG_{\pi} = \{\G_{\ell} = p_1p_2 \cdots p_{\ell} \mZ \mid \ell \geq 1\}$. Let $\whmZ_{\pi}$ be the completion of $\mZ$ with respect to the chain $\cG_{\pi}$. Then we have a topological isomorphism
 \begin{equation}
\whmZ_{\pi} \cong \prod_{i \geq 1} \ \mZ/p_i \mZ \ .
\end{equation}
Let $H_{\ell} = \mZ/p_1\mZ \oplus \cdots \oplus \mZ/p_{\ell} \mZ$ be the direct sum of the first $\ell$-factors. Then $\{H_{\ell} \mid \ell \geq 1\}$ is an infinite increasing chain of finite subgroups of $\whmZ_{\pi}$   which does not stabilize, so $\whmZ_{\pi}$ is not topologically Noetherian.

  These two examples illustrate  the   idea behind the proof of the following result. 
\begin{prop}\label{prop-nilpNoetherian}
Let $\G$ be a finitely generated nilpotent group, and let $\whGamma$ be a profinite completion of $\G$.
Then $\whGamma$ is topologically Noetherian if and only if   the prime spectrum $\pi(\Pi[\whGamma])$ is finite.
\end{prop}
\proof
Recall some basic facts about profinite groups. (See for example,  \cite[Chapter~2]{Wilson1998}.) For a prime $p$,  a finite group $H$ is a $p$-group if every element of $H$ has order a power of $p$. A profinite group $\fH$ is a pro-$p$-group if $\fH$ is the inverse limit of finite $p$-groups. A Sylow $p$-subgroup $\fH \subset \fG$ is a maximal pro-$p$-subgroup \cite[Definition~2.2.1]{Wilson1998}. 

If $\fG$ is pro-nilpotent, then for each prime $p$, there is a unique Sylow $p$-subgroup of $\fG$, which is normal in $\fG$ \cite[Proposition~2.4.3]{Wilson1998}. Denote this group by $\fG_{(p)}$. Moreover, $\fG_{(p)}$ is non-trivial if and only if 
$p \in \pi(\Pi[\fG])$. It follows that there 
 is a topological isomorphism  
\begin{equation}\label{eq-primeSylowdecomp}
\fG \cong \prod_{p \in \pi(\Pi[\fG])} ~ \fG_{(p)} \ .
\end{equation}
 From the isomorphism \eqref{eq-primeSylowdecomp} it follows immediately that if the prime spectrum $\pi(\Pi[\fG])$ is infinite, then $\fG$ is not topologically Noetherian. To see this, list $\pi(\Pi[\fG]) = \{p_i \mid i = 1,2, \ldots \}$, then we obtain an infinite strictly  increasing chain of closed subgroups,
 $$\fH_{\ell} = \prod_{i=1}^{\ell} \ \fG_{(p_i)} \ . $$
 If the prime spectrum $\pi(\Pi[\fG])$ is finite,  then the isomorphism   \eqref{eq-primeSylowdecomp} reduces the proof that $\fG$ is topologically Noetherian to the case of showing that if $\fG$ is topologically finitely generated, then each of its Sylow $p$-subgroups is Noetherian. The group  $\fG_{(p)}$ is nilpotent and topologically finitely generated, so we can use  the lower central series for $\fG_{(p)}$ and induction to reduce to the case where 
 $\fH$ is a topologically finitely-generated abelian pro-$p$-group, and so is isomorphic to a finite product of $p$-completions of $\mZ$, which are topologically Noetherian.
  
The proof of Proposition~\ref{prop-nilpNoetherian} is completed by observing that   a   profinite completion $\whGamma$ of a finitely generated nilpotent group $\G$ is a topologically finitely-generated nilpotent group, and we apply the above remarks.
\endproof

\begin{cor}\label{cor-vnilpotentNoetherian}
Let $\G$ be a virtually nilpotent group, that is there exists a finitely-generated nilpotent subgroup $\G_0 \subset \G$ of finite index.
Then   a    profinite completion $\whGamma$ of $\G$  is topologically Noetherian if and only if   its prime spectrum $\pi(\Pi[\whGamma])$ is finite.
\end{cor}
\proof
We can assume that $\G_0$ is a normal subgroup of $\G$, then its closure $\whGamma_0 \subset \whGamma$ satisfies the hypotheses of Proposition~\ref{prop-nilpNoetherian}, and the Steinitz orders satisfy $[\whGamma_0] \mor [\whGamma]$.
As $\whGamma_0$ is topologically Noetherian if and only if $\whGamma$ is topologically Noetherian, the claim follows.
\endproof

\subsection{Dynamics of Noetherian groups}\label{subsec-Noetheriandynamics}
We next relate the topologically Noetherian property of a profinite group with the dynamics of a Cantor action of the group, to obtain proofs of Theorem~\ref{thm-nilstable} and Corollary~\ref{cor-niltopfree}.  We first give  the profinite analog of     \cite[Theorem~1.6]{HL2018b}. We follow the  outline of its proof.
  
\begin{prop}\label{prop-NLQA}
Let $\fG$ be a topologically Noetherian group.
Then     a minimal equicontinuous action   $(\fX,\fG,\whPhi)$  on a Cantor space $\fX$   is locally quasi-analytic.
\end{prop}
\proof
First, the closure $\fG(\Phi)  \subset \Homeo(\fX)$,   so     the action $\whPhi$ of $\fG(\Phi)$ is effective.  Suppose that the action $\whPhi$ is not locally quasi-analytic, then there exists an infinite properly decreasing chain of clopen subsets of $\fX$,
$\{U_1 \supset U_2 \supset \cdots \}$, which satisfy the following properties, for all $\ell \geq 1$:
\begin{itemize}
\item $U_{\ell}$ is adapted to the action $\whPhi$ with isotropy subgroup $\fG_{U_{\ell}} \subset \fG$;
\item there is a closed subgroup $K_{\ell} \subset \fG_{U_{\ell+1}}$ whose restricted action to $U_{\ell +1}$ is trivial, but the restricted action of $K_{\ell}$ to $U_{\ell}$ is effective.
\end{itemize}
It follows that we obtain a properly increasing chain of closed subgroups $\{K_1 \subset K_2 \subset \cdots\}$ in $\fG$, which contradicts the assumption that $\fG$ is topologically Noetherian.
\endproof

We now give the  proof  of Theorem~\ref{thm-nilstable}.
Let $(\fX,\G,\Phi)$ be a  nilpotent Cantor action. 
 Then there exists a finitely-generated  nilpotent   subgroup $\G_0 \subset \G$ of finite index, and we can assume without loss  of generality   that $\G_0$ is normal. Let $\whGamma_0$ be the closure of $\G_0$ in $\whGamma$ and let $x \in \fX$ be a basepoint. Note that the group $\whGamma$ has finite prime spectrum   if and only if the group $\whGamma_0$ has finite prime spectrum. Thus,  it suffices to show that the action of $\G_0$ on the orbit $\fX_0 = \whGamma_0 \cdot x$ is stable. For simplicity of notation, we will simply assume that the given group $\G$ is itself nilpotent.

The profinite completion $\fG(\G)$  of $\Phi(\G)$ is also nilpotent, and we have the profinite action $(\fX,\whGamma,\whPhi)$.
 Suppose that the action $\whPhi$ is not stable, then   there exists an increasing  chain of closed subgroups $\{K_{\ell} \mid \ell \geq 1\}$ where $K_{\ell}$ acts trivially on the clopen subset $U_{\ell} \subset \fX$. Let $x \in \cap_{\ell > 0} \ U_{\ell}$ then each $K_{\ell} \subset \fD(\Phi, x)$, so $\fD(\Phi, x)$ contains a strictly increasing chain of closed subgroups. As we are given that  the prime spectrum $\pi(\Pi[\fD(\Phi, x)])$ is     finite, this contradicts the 
  conclusion  of  Proposition~\ref{prop-nilpNoetherian}. Hence, the action $\whPhi$ must be locally quasi-analytic, as was to be shown.

The proof of   Corollary~\ref{cor-niltopfree} is just an extension of that of Theorem~\ref{thm-nilstable}.
 Let $(\fX,\G,\Phi)$ be a nilpotent Cantor action for which the Steinitz order $\Pi(\fG(\Phi))$ has prime multiplicities at most $2$, at all but a finite number of primes. 
 As before, we can assume without loss of generality  that the group $\G$ is nilpotent. Then we have the decomposition \eqref{eq-primeSylowdecomp} of $\fG(\Phi)$ into a product of its Sylow $p$-subgroups, and the corresponding product decomposition of the space
 \begin{equation}\label{eq-Xproduct}
\fX \cong  \prod_{p \in \pi(\Pi[\fG(\Phi)])} ~ \fX_{(p)} ~ = ~ \prod_{p \in \pi(\Pi[\fG(\Phi)])} ~ \fG(\Phi)_{(p)}/\fD(\Phi)_{(p)} \ .
\end{equation}
The factors in the product  representation of $\fG(\Phi)$ in  \eqref{eq-primeSylowdecomp} act on the corresponding factors in \eqref{eq-Xproduct}.
In particular, 
 the factors $\fG(\Phi)_{(p)}$ and $\fG(\Phi)_{(q)}$ commute when $p \ne q$, and thus their actions on $\fX$ commute.
Also note that if  the multiplicity of $p$ is finite, then the corresponding Sylow $p$-subgroup $\fG(\Phi)_{(p)}$ is a finite group, and so the quotient space $\fX_{(p)}$ is a finite set.
 
 Let   $\fG(\Phi)_{(p)}$ be a $p$-Sylow subgroup with order at most $p^2$. Then $\fG(\Phi)_{(p)}$ is a nilpotent group of order $p^2$, so must   be abelian.

Let $\cD(\Phi)$ denote the discriminant of the action $\Phi$. Its  $p$-Sylow subgroup satisfies $\cD(\Phi)_{(p)} \subset \fG(\Phi)_{(p)}$.  

 If the multiplicity of $p$ is at most $2$, then for $\whg \in  \cD(\Phi)$, the left action of its projection to $\cD(\Phi)_{(p)}$ fixes the basepoint in $\fX_{(p)}$, and as $\cD(\Phi)_{(p)}$ is abelian, the action fixes 
all of the points in  the finite quotient space $\fX_{(p)} = \fG(\Phi)_{(p)}/\fD(\Phi)_{(p)}$. As the action of a non-trivial element of $\fD(\Phi)_{(p)}$ must be non-trivial, this implies the projection is the identity element in $\fG(\Phi)_{(p)}$.

Thus, it suffices to show that the action of $\whg$ on the factors in \eqref{eq-Xproduct}
for which the prime order $n(p) \geq 3$ is stable. As there are at most a finite number of such factors, we are reduced to the situation in the proof of Theorem~\ref{thm-nilstable}, and so the action must be stable.

 \section{Examples}\label{sec-examples}

We give in this section a collection of examples of nilpotent Cantor actions to illustrate the results and ideas  of this work. Our guiding principle is to present the simplest examples in each class, which can then be made as complicated as desired following the basic design. All of these examples give rise to  solenoidal manifolds with the specified prime spectrum, with base manifold   an $n$-torus in Example~\ref{ex-toral}, or base manifold the standard compact nil-3-manifold for Examples~\ref{ex-trivial}, \ref{ex-stable} and \ref{ex-wild}.

 \subsection{Toroidal actions}\label{subsec-toral}
 
We begin with the simplest examples of Cantor actions for which the prime spectra are not sufficient to distinguish the actions.
  A \emph{toroidal Cantor action} is the action  of $\G = \mZ^m$ on a ``diagonal'' profinite completion of $\mZ^m$, for some $m \geq 1$. The classification of minimal equicontinuous actions of $\mZ^m$ involves subtleties associated with the space of lattice chains in $\mR^m$, as discussed in various works \cite{GPS2019,Li2018}. The diagonal actions, which we now define, suffice for illustrating the construction of actions with prescribed prime spectrum. 

 \begin{ex}\label{ex-toral}
 {\rm 
Consider the case $n=1$. Choose two disjoint sets of distinct primes,
$$\pi_f = \{q_1 , q_2, \ldots \} \quad , \quad \pi_{\infty} = \{p_1 , p_2, \ldots\}$$
where $\pi_f$ and $\pi_{\infty}$ can be chosen to be finite or infinite sets, and either $\pi_f$ is infinite, or $\pi_{\infty}$ is non-empty. Choose  multiplicities $n(q_i) \geq 1$ for the primes in $\pi_f$.
For each $\ell > 0$, define a subgroup of $\G = \mZ$ by 
$$\G_{\ell} = \{q_1^{n(q_1)} q_2^{n(q_2)} \cdots q_{\ell}^{n(q_{\ell})} \cdot p_1^{\ell} p_2^{\ell} \cdots p_{\ell}^{\ell} \cdot n \mid n \in \mZ \} \ , $$
with the understanding that if the prime $q_{\ell}$ or $p_{\ell}$ is not defined, then we simply set this term to be $1$.
The completion $\whGamma$ of $\mZ$ with respect to this group chain  admits a product decomposition into its Sylow $p$-subgroups
\begin{equation}\label{eq-pqlimit}
\whGamma ~ \cong ~ \prod_{i =1}^{\infty} \  \mZ/q_i^{n(q_i)} \mZ ~   \cdot ~ \prod_{p \in \pi_{\infty}} ~  \whmZ_{(p)} \ ,
\end{equation}
where $\whmZ_{(p)}$ denotes the $p$-adic completion of $\mZ$.
Thus   $\pi(\Pi[\whGamma]) = \pi_f \cup \pi_{\infty}$. As $\mZ$ is abelian,   we have $\fX = \whGamma$ and the the discriminant group for the action of $\G$   is trivial.
}
\end{ex}

 \begin{ex}\label{ex-almosttoral}
 {\rm 
We next give two extensions of the diagonal actions described in Example~\ref{ex-toral}. 

First, we construct a 
 diagonal toroidal action of $\mZ^m$   by making $m$ choices of prime spectra as above, then taking the product action. While the return equivalence class of a $\mZ$-action on $\fX = \whGamma$ as in \eqref{eq-pqlimit} is determined by the asymptotic class $\Pi_a[\whGamma]$, as in Theorem~\ref{thm-onedimSol}, this need no longer hold  for the product of such actions. 
For example, the two profinite completions of $\mZ^2 = \mZ \oplus \mZ$ given by
\begin{equation}
\whGamma_1 =  \whmZ_{(6)} \oplus \whmZ_{(5)}  \quad , \quad \whGamma_2 = \whmZ_{(2)} \oplus \whmZ_{(15)}
\end{equation}
have the same Steinitz orders, but are not isomorphic.

The second construction shows that the conclusion of Theorem~\ref{thm-returnequivspectra} is best possible,   that is, return equivalence need not preserve the Steinitz order of the action. 
Let $\pi_f = \{p_1, p_2, \ldots \}$ be a proper subset of  primes, infinite in number and all distinct. Let $\whmZ_{\pi_f}$ denote the completion of $\mZ$ with respect to the primes $\pi_f$ where we choose multiplicity $n(p) = 1$ for each $p \in \pi_f$.
Then we have the odometer action $\Phi_1$ of $\mZ$ on $\fX_1 = \whmZ_{\pi_f}$.

Next, for $k \geq 2$,  consider the action of $\mZ^k = \mZ \oplus \cdots \oplus \mZ$ on $\fX = \whmZ_{\pi_f} \oplus \cdots \oplus \whmZ_{\pi_f}$.
Let $\G = \mZ^k \rtimes C_k$  where $C_k = \mZ/k \mZ$ is the cyclic group of order $k$, which acts on the   factor $\mZ^k$  by the automorphism which is a cyclic permutation of the basis vectors. Then $C_k$ also acts on $\fX$ by the corresponding cyclic permutation of the factors, and we use this to define an action $\Phi_2$ of $\G$ on $\fX$.

The actions $\Phi_1$ and $\Phi_2$ are return equivalent. To see this, observe that the coset of the identity in $C_k$ determines a clopen subset of $\fX$, and the restriction of the action $\Phi_2$ to this coset is just the odometer action $\Phi_1$.

Suppose that $k$ is a prime which is not in $\pi_f$, then $\pi(\Pi[\Phi_2]) = \pi_f \cup  \{k\} = \pi(\Pi[\Phi_1]) \cup \{k\}$,   and so their prime spectra differ. If $k$ is a prime which is   in $\pi_f$ then the prime spectra of the two actions agree, but their multiplicities do not.
 One  can also repeat this construction for any transitive subgroup of the permutation group $\Perm(k)$ on $k$-elements for $k \geq 2$, and so obtain that  the prime spectra of the two actions differ by an arbitrary set of primes  which are divisors of $k$. 
}
 \end{ex}

 \subsection{Heisenberg actions}\label{subsec-heisenberg}

  We next construct a selection of examples, given by the action of the integer Heisenberg group $\cH$ on a profinite completion of the group. The group $\cH$ is a cocompact lattice in the real Heisenberg group $H_3(\mR)$, so the quotient $M = H_3(\mR)/\cH$ is a compact $3$-manifold, and the choice of a group chain in $\cH$ defines a tower of coverings of $M$ whose inverse limit has monodromy action conjugate to the Cantor actions defined by the group chain.

Let $\cH$ be  represented as the upper triangular matrices in ${\rm GL(\mZ^3)}$. That is,
\begin{equation}\label{eq-cH}
\cH =   \left\{  \left[ {\begin{array}{ccc}
   1 & a & c\\
   0 & 1 & b\\
  0 & 0 & 1\\
  \end{array} } \right] \mid a,b,c  \in \mZ\right\} .
\end{equation}
In coordinates $(a,b,c) , (a',b',c') \in \mZ^3$, the group operation  $*$ and inverse  are given by, 
\begin{equation}\label{eq-Hrules}
(a,b,c)*(a',b',c')=(a+a',b+b',c+c'+ab') \quad , \quad (a,b,c)^{-1} = (-a, -b, -c +ab) \ .
\end{equation}
  In particular, we have 
\begin{equation}\label{eq-Hcomm}
(a,b,c) * (a',b',c')*(a,b,c)^{-1}=(a',b',c' + ab' -ba') \ .
\end{equation}
The work \cite{LSU2014} gives a complete discussion of the normal subgroups of $\cH$.

 \begin{ex}\label{ex-trivial}
 {\rm
 We construct a Cantor action of $\cH$ on a profinite completion defined by   a proper self-embedding of $\cH$ into  itself. The resulting action has trivial discriminant group, but the integers $k_{\ell}$ and $k^*_{\ell}$ defined in \eqref{eq-dims} are distinct.
 The variety of such actions have been extensively studied in the authors' work \cite{HLvL2020}  joint with Van Limbeek, as they all yield stable Cantor actions.

For a  prime $p \geq 2$, define the self-embedding $\vp_p \colon \cH \to \cH$ by  
$\vp(a,b,c) = (pa, pb, p^2c)$. Then define a group chain in $\cH$ by setting 
$$\cH_{\ell} = \vp_p^{\ell}(\cH) = \{(p^{\ell} a, p^{\ell}b, p^{2\ell}c) \mid a,b,c \in \mZ\} \quad, \quad \bigcap_{\ell > 0} \ \cH_{\ell} = \{e\} \ .$$

Formula \eqref{eq-Hcomm} implies that the normal core for $\cH_{\ell}$ is given by
$$C_{\ell} = {\rm core}(\cH_{\ell})  = \{(p^{2\ell} a, p^{2\ell} b, p^{2\ell} c) \mid a,b,c \in \mZ\} \ .$$
Thus,  the finite group 
\begin{equation}\label{eq-Qell}
Q_{\ell} = \cH/C_{\ell} \cong \{( \oa, \ob, \oc) \mid \oa, \ob, \oc\in \mZ/p^{2\ell}\mZ \} \ .
\end{equation}
The profinite group $\widehat{\cH}_{\infty}$ is the inverse limit of the quotient groups $Q_{\ell}$ so we have
$$\widehat{\cH}_{\infty} =   \{(\wha, \whb, \whc) \mid \wha, \whb, \whc\in \widehat{\mZ}_{p^2} \}$$
 with multiplication  on each finite quotient  induced  by the formula \eqref{eq-Hrules}.
 Note that the group $\cH$ embeds into $\widehat{\cH}_{\infty}$ as $p^{\ell}$ tends to infinity with $\ell$.
 
Next, we calculate the discriminant subgroup $\cD_{\infty}$ for this action. First note   
\begin{eqnarray}
\cH_{\ell}/C_{\ell} & = & \{(p^{\ell} \oa, p^{\ell} \ob, 0) \mid \oa, \ob \ \in \mZ/p^{\ell}\mZ \} \ \subset \ Q_{\ell} \ ,  \label{eq-Qell1}\\
\cH_{\ell+1}/C_{\ell+1} & = & \{(p^{\ell +1} \oa, p^{\ell +1} \ob, 0) \mid \oa, \ob \ \in   \mZ/p^{\ell +1}\mZ \} \   \ . \label{eq-Qell=2}
\end{eqnarray}
 Thus, $k_{\ell} = \# (\cH_{\ell}/C_{\ell}) = p^{2\ell}$.

  Note that  $\cH_{2\ell} \subset C_{\ell}$. So while each quotient $\cH_{2\ell}/C_{2\ell}$ is non-trivial, its image under the composition of bonding maps in 
   \eqref{eq-discformula} vanishes in $\cH_{\ell}/C_{\ell}$.  Thus $\cD_{\vp}$ is the trivial group, and so  each $k_{\ell}^* = 1$.   

 }
 \end{ex}

\begin{ex}[A toy model] \label{ex-toy} 
{\rm
We describe a  \emph{finite}  action which is used to construct the next classes of Heisenberg actions which have non-trivial discriminant groups, and arbitrary prime spectra.  
   
   Fix a prime $p \geq 2$. For $n \geq 1$ and  $0 \leq k < n$, we have the following finite groups: 
 \begin{equation}
G_{p,n} =   \left\{  \left[ {\begin{array}{ccc}
   1 & \oa & \oc\\
   0 & 1 & \ob\\
  0 & 0 & 1\\
  \end{array} } \right] \mid \oa,\ob,\oc  \in \mZ/p^{n}\mZ\right\} ~ , ~
  H_{p,n,k} =   \left\{  \left[ {\begin{array}{ccc}
   1 & p^k \oa & 0\\
   0 & 1 & 0\\
  0 & 0 & 1\\
  \end{array} } \right] \mid \oa  \in \mZ/p^{n}\mZ\right\}
\end{equation}
Note that $\#[G_{p,n}] = p^{3n}$ and $\#[H_{p,n,k}] = p^{n-k}$.

Let   $\ox = (1,0,0), \oy = (0,1,0), \oz= (0,0,1) \in G_{p,n}$, then by formula \eqref{eq-Hcomm} we have
$\ox \cdot \oy \cdot \ox^{-1} = \oy \oz$ and $\ox \cdot \oz \cdot \ox^{-1} = \oz$.  That is, the adjoint action of $\ox$ on the ``plane'' in the  $(\oy,\oz)$-coordinates   is a ``shear'' action along the $\oz$-axis, and the adjoint action of $\ox$ on the $\oz$-axis fixes   all points on the $\oz$-axis.

Set $X_{p,n,k} = G_{p,n}/H_{p,n,k}$,  then the isotropy group of the action of $G_{p,n}$ on $X_{p,n,k}$ at the coset  $H_{p,n,k}$  of the identity element is $H_{p,n,k}$. The core subgroup $C_{p,n,k} \subset H_{p,n,k}$ contains elements in $H_{p,n,k}$ which fix every point in $X_{p,n,k}$. The action of $\ox$ on the coset space $X_{p,n,k}$ satisfies $\Phi(\ox)(\oy) = \oy \oz$, so the identity is the only element in $H_{p,n,k}$, so $C_{p,n,k}$ is trivial. Then $D_{p,n,k} = H_{p,n,k}/C_{p,n,k} = H_{p,n,k}$, and for each $g \in H_{p,n,k}$ its action fixes the multiples of $\oz$. 
}
\end{ex}

In the following two classes of examples,    given sets of primes $\pi_f$ and $\pi_{\infty}$, we embed a   infinite product  of finite actions as in Example~\ref{ex-toy} into   a profinite completion $\whcH_{\infty}$ of $\cH$, which defines a nilpotent Cantor action $(X_{\infty}, \cH, \Phi)$  on a quotient $X_{\infty} = \whcH_{\infty}/D_{\infty}$.  
This is possible, due to the following result for pro-nilpotent groups, which is a consequence of \cite[Proposition~2.4.3]{Wilson1998}. 
\begin{prop}\label{prop-factorization}
Let $\whGamma$ be a profinite completion of a finitely-generated nilpotent group $\G$. Then there is a topological isomorphism
\begin{equation}\label{eq-factorization}
\whGamma \cong \prod_{p \in \pi(\Pi[\whGamma])} \ \whGamma_{(p)} \ ,
\end{equation}
where $\whGamma_{(p)} \subset \whGamma$ denotes the Sylow $p$-subgroup of $\whGamma$ for a prime $p$.
\end{prop}

\begin{ex}[Stable Heisenberg actions] \label{ex-stable} 
{\rm
 We   construct     Heisenberg actions  with finite or infinite prime spectrum, using the product formula \eqref{eq-factorization}, and then show that they are stable. 
 
Let $\pi_f $ and  $\pi_{\infty}$ be two disjoint   collections of primes, with $\pi_f$ a finite set, and $\pi_{\infty}$ a non-empty set.   
Enumerate $\pi_f = \{q_1, q_2, \ldots, q_m\}$ then choose integers $1 \leq r_i \leq n_i$ for $1 \leq i \leq m$. Enumerate $\pi_{\infty} = \{p_1, p_2, \ldots\}$  with the convention (for notational convenience) that if $\ell$ is greater than the number of primes in $\pi_{\infty}$ then we set $p_{\ell} = 1$.
For each $\ell \geq 1$,  define the integers
\begin{eqnarray}
M_{\ell}    & =   &  q_1^{r_1} q_2^{r_2} \cdots q_m^{r_m} \cdot p_1^{\ell} p_2^{\ell} \cdots p_{\ell}^{\ell} \  , \\
N_{\ell}    & =   &  q_1^{n_1} q_2^{n_2} \cdots q_m^{n_m} \cdot p_1^{\ell} p_2^{\ell} \cdots p_{\ell}^{\ell} \  ,
\end{eqnarray}
 
For all $\ell \geq 1$, observe that $M_{\ell}$ divides $N_{\ell}$, and    define a subgroup of $\cH$, in the coordinates above, 
\begin{equation}\label{eq-chain}
\cH_{\ell} = \{ (a M_{\ell},b N_{\ell} ,c N_{\ell}) \mid a,b,c \in \mZ \} \ .
\end{equation}
Its  core  subgroup is given by 
$C_{\ell} =  \{ (a N_{\ell},b N_{\ell} ,c N_{\ell}) \mid a,b,c \in \mZ \}$. 
Observe that
$$\mZ/N_{\ell} \mZ \cong \mZ/q_1^{n_1}\mZ \oplus \cdots \oplus  \mZ/q_m^{n_m}\mZ \oplus \mZ/p_1^{\ell}\mZ \oplus \cdots \oplus \mZ/p_{\ell}^{\ell}\mZ \ .$$
 By Proposition~\ref{prop-factorization},  and in the notation of  Example~\ref{ex-toy},  we have for   $k_i = n_i - r_i$ that
\begin{equation}\label{eq-lqafactors}
 \whcH_{\infty} ~ \cong ~ \prod_{i=1}^m \ G_{q_i, n_i} ~ \cdot ~ \prod_{j=1}^{\infty} \ \whcH_{(p_j)} \quad , \quad D_{\infty} ~ \cong ~  \prod_{i=1}^m \ H_{q_i, n_i, k_i} \ . 
\end{equation}
Then the Cantor space  $X_{\infty} = \whcH_{\infty}/D_{\infty}$ associated to the group chain $\{\cH_{\ell} \mid \ell \geq 1\}$ is given by
\begin{equation}\label{eq-lqaspace}
X_{\infty} ~ \cong ~ \prod_{i=1}^m \ X_{q_i, n_i, k_i} ~ \times ~ \prod_{j=1}^{\infty} \ \whcH_{(p_j)}  \ . 
\end{equation}

In particular,  as the first factor in   \eqref{eq-lqaspace} is a finite product of finite sets, the second factor defines an open neighborhood $$U = \prod_{i=1}^m \ \{x_i\} ~ \times ~ \prod_{j=1}^{\infty} \ \whcH_{(p_j)}$$ 
where $x_i \in X_{q_i, n_i, k_i}$ is the basepoint given by the coset of the identity element.
That is, $U$ is a clopen neighborhood of the basepoint in $X_{\infty}$.  The isotropy group of $U$ is given by 
\begin{equation}
 \whcH_{\infty}|U ~ = ~  \prod_{i=1}^m \ H_{q_i, n_i, k_i} ~ \times ~ \prod_{j=1}^{\infty} \ \whcH_{(p_j)}  \ . 
\end{equation}
The restriction of  $\whcH_{\infty}|U$ to $U$ is isomorphic to the subgroup  
\begin{equation}
 K|U ~ = ~  \prod_{i=1}^m \ \{\overline{e}_i\} ~ \times ~ \prod_{j=1}^{\infty} \ \whcH_{(p_j)} ~ \subset ~ \Homeo(U) \ , 
\end{equation}
where $\overline{e}_i \in G_{q_i, n_i}$ is the identity element.  The group $K|U$   acts freely on $U$, and thus the action of $\whcH_{\infty}$ on $X_{\infty}$ is locally quasi-analytic.
Moreover, the union $\pi = \pi_f \cup \pi_{\infty} = \pi(\Pi[\whcH_{\infty}])$ is the prime spectrum of the action of $\cH$ on $X_{\infty}$.  If $\pi_\infty$ is infinite, then the prime spectrum of the action is infinite.
Note that the group $\cH$ embeds into $\widehat{\cH}_{\infty}$ as the integers $M_{\ell}$ and $N_{\ell}$ tend to infinity with $\ell$.

 }
\end{ex}

\begin{ex}[Wild Heisenberg actions]\label{ex-wild} 
{\rm

Let $\pi_f $ and  $\pi_{\infty}$ be two disjoint   collections of primes, with $\pi_f$ an infinite set and $\pi_{\infty}$ arbitrary, possibly empty.   
Enumerate $\pi_f = \{q_1, q_2, \ldots\}$ and choose integers $1 \leq r_i \leq n_i$ for $1 \leq i < \infty$. Enumerate $\pi_{\infty} = \{p_1, p_2, \ldots\}$, again with the convention   that if $\ell$ is greater than the number of primes in $\pi_{\infty}$ then we set $p_{\ell} = 1$. 

As in Example~\ref{ex-stable}, for  each $\ell \geq 1$,  define the integers
\begin{eqnarray}
M_{\ell}    & =   &  q_1^{r_1} q_2^{r_2} \cdots q_{\ell}^{r_{\ell}} \cdot p_1^{\ell} p_2^{\ell} \cdots p_{\ell}^{\ell} \  , \\
N_{\ell}    & =   &  q_1^{n_1} q_2^{n_2} \cdots q_{\ell}^{n_{\ell}} \cdot p_1^{\ell} p_2^{\ell} \cdots p_{\ell}^{\ell} \  .
\end{eqnarray}

For $\ell \geq 1$, define a subgroup of $\cH$, in the coordinates above, 
\begin{equation}\label{eq-chain2}
\cH_{\ell} = \{ (a M_{\ell},b N_{\ell} ,c N_{\ell}) \mid a,b,c \in \mZ \} \ , 
\end{equation}
 Its core  subgroup is given by 
$C_{\ell} =  \{ (a N_{\ell},b N_{\ell} ,c N_{\ell}) \mid a,b,c \in \mZ \}$. For   $k_i = n_i - r_i$ we then have
\begin{equation}\label{eq-lqafactors2}
 \whcH_{\infty} ~ \cong ~ \prod_{i=1}^{\infty} \ G_{q_i, n_i} ~ \cdot ~ \prod_{j=1}^{\infty} \ \whcH_{(p_j)} \quad , \quad D_{\infty} ~ \cong ~  \prod_{i=1}^{\infty} \ H_{q_i, n_i, k_i} \ . 
\end{equation}
The Cantor space  $X_{\infty} = \whcH_{\infty}/D_{\infty}$ associated to the group chain $\{\cH_{\ell} \mid \ell \geq 1\}$ is given by
\begin{equation}\label{eq-lqaspace2}
X_{\infty} ~ \cong ~ \prod_{i=1}^{\infty} \ X_{q_i, n_i, k_i} ~ \times ~ \prod_{j=1}^{\infty} \ \whcH_{(p_j)}  \ . 
\end{equation}
The first factor in   \eqref{eq-lqaspace} is an infinite product of finite sets, so fixing the first $\ell$-coordinates in this product determines a clopen subset of $X_{\infty}$. Let $x_i \in X_{q_i, n_i, k_i}$ denote the  coset of the identity element, which is   the basepoint in $X_{q_i, n_i, k_i}$. Then for  each $\ell \geq 1$, we define a clopen set in $X_{\infty}$ 
\begin{equation}\label{eq-wildclopen}
U_{\ell} = \prod_{i=1}^{\ell} \ \{x_i\} ~ \times ~  \prod_{i=\ell+1}^{\infty} \ X_{q_i, n_i, k_i} ~ \times ~ \prod_{j=1}^{\infty} \ \whcH_{(p_j)} \ .
\end{equation}
   Recalling the calculations in   Example~\ref{ex-toy}, the subgroup $H_{q_i, n_i, k_i}$ is the isotropy group of the basepoint $x_i \in X_{q_i, n_i, k_i}$. Thus, 
the isotropy subgroup of $U_{\ell}$  for the $\whcH_{\infty}$-action  is given by the product
\begin{equation}
 \whcH_{\infty}|_{U_{\ell}} ~ = ~  \prod_{i=1}^{\ell} \ H_{q_i, n_i, k_i} ~ \times ~  \prod_{i=\ell + 1}^{\infty} \ G_{q_i, n_i} ~ \times ~ \prod_{j=1}^{\infty} \ \whcH_{(p_j)}  \ . 
\end{equation}
For $j \ne i$, the subgroup $H_{q_i, n_i, k_i}$ acts as the identity on the factors $X_{q_j, n_j, k_j}$ in \eqref{eq-lqaspace2}.
Thus, the  image of  $\whcH_{\infty}|_{U_{\ell}}$ in $\Homeo(U_{\ell})$  is isomorphic to the subgroup  
\begin{equation}
 Z_{\ell} ~ = ~\whcH_{\infty}|U_{\ell} ~ = ~  \prod_{i=1}^{\ell} \ \{\overline{e}_i\} ~ \times ~  \prod_{i=\ell + 1}^{\infty} \ G_{q_i, n_i} ~ \times ~ \prod_{j=1}^{\infty} \ \whcH_{(p_j)}  ~ \subset ~ \Homeo(U_{\ell}) \ , 
\end{equation}
where $\overline{e}_i \in G_{q_i, n_i}$ is the identity element.  

We next show that this action   is not stable; that is, for any $\ell > 0$ there exists a clopen subset $V \subset U_{\ell}$ and non-trivial $\whg \in   Z_{\ell}$ so that the action of $\whG$ restricts to the identity map on $V$. We can assume without loss that $V= U_{\ell'}$  for some $\ell' > \ell$. Consider  the restriction map for 
the isotropy subgroup of $Z_{\ell}$ to  $U_{\ell'}$ which is given by 
  $$\rho_{\ell, \ell'} \colon Z_{\ell}|_{U_{\ell'}}  \to Z_{\ell'} \subset \Homeo(U_{\ell'}) \ .$$
We must show that there exists $\ell' > \ell$ such that this map has    a non-trivial kernel.  
 Calculate this map in terms of the product representations above,
  \begin{equation}\label{eq-restriction}
 Z_{\ell}|_{U_{\ell'}}  ~ = ~  \prod_{i=1}^{\ell} \ \{\overline{e}_i\} ~ \times ~  \prod_{i=\ell + 1}^{\ell'} \ H_{q_i, n_i,k_i} ~ \times ~  \prod_{i=\ell' + 1}^{\infty} \ G_{q_i, n_i}~ \times ~ \prod_{j=1}^{\infty} \ \whcH_{(p_j)}   \ . 
\end{equation}
For $\ell < i \leq \ell'$, the group $H_{q_i, n_i,k_i}$ fixes the point $\prod_{i=1}^{\ell'} \ \{x_i\}$, and acts trivially on 
$\prod_{i=\ell'+1}^{\infty} \ X_{q_i, n_i, k_i}$.
Thus, the kernel of the restriction map contains the second factor in \eqref{eq-restriction}, 
\begin{equation}
 \prod_{i=\ell + 1}^{\ell'} \ H_{q_i, n_i,k_i} ~ \subset ~  \ker \left\{ \rho_{\ell, \ell'} \colon Z_{\ell}|_{U_{\ell'}}  \to \Homeo(U_{\ell'})  \right\} \ .
\end{equation}
As this group is non-trivial for all $\ell' > \ell$,  the action of $\whcH_{\infty}$ on $X_{\infty}$ is not locally quasi-analytic, hence the action of $\cH$ on $X_{\infty}$ is wild.
Moreover, the prime spectrum of the action of $\cH$ on $X_{\infty}$ equals the union $\pi = \pi_f \cup \pi_{\infty}$.
 }
 \end{ex}

Finally, we give the proof of Theorem~\ref{thm-lqaWILD} using the construction in Example~\ref{ex-wild},  that is, we show that choices in Example \ref{ex-wild} can be made in such a way that the action of $\cH$ on a Cantor set is topologically free while the action of $\whcH_\infty$ is not stable. To do that, choose  an infinite set of distinct primes $\pi_f = \{q_1, q_2, \ldots\}$, and let the set of infinite primes $\pi_{\infty}$ be empty.
Choose the constants $n_i =2$ and $k_i=1$ for all $i \geq 1$. Let $X_{\infty}$ be the Cantor space  defined by \eqref{eq-lqaspace2}. Then the action of $\cH$  is   wild by the calculations in Example~\ref{ex-wild}. 

We claim that this action is   topologically free. Suppose not, then there exists an open set $U \subset X_{\infty}$ and $g \in \cH$ such that the action of $\Phi(g)$ is non-trivial on   $X_{\infty}$ but leaves the set $U$ invariant, and restricts to the identity action on $U$.
The action of $\cH$ on $X_{\infty}$ is minimal, so there exists $h \in \cH$ with $h \cdot x_{\infty} \in U$. Then $\Phi(h^{-1} g h)(x_{\infty}) = x_{\infty}$ and the action   $\Phi(h^{-1} g h)$ fixes an open neighborhood of $x_{\infty}$. Replacing $g$ with $h^{-1} g h$ we can assume that $\Phi(g)(x_{\infty}) = x_{\infty} \in U$. From the definition \eqref{eq-wildclopen}, the clopen sets
\begin{equation}\label{eq-wildclopen2}
U_{\ell} = \prod_{i=1}^{\ell} \ \{x_i\} ~ \times ~  \prod_{i=\ell+1}^{\infty} \ X_{q_i, 2, 1}  
\end{equation}
form a neighborhood basis at $x_{\infty}$, and thus there exists $\ell > 0$ such that $U_{\ell} \subset U$.

The group $\cH$  diagonally embeds into $\whcH_{\infty}$ so from the expression \eqref{eq-lqafactors2}, we have
 $\ds g = (g,g,\ldots,g) \in    \prod_{i=1}^{\infty} \ G_{q_i, 2}$. The action of $\Phi(g)$ is factorwise, and $\Phi(g)(x_{\infty}) = x_{\infty}$ implies that $\ds g \in D_{\infty} \cong   \prod_{i=1}^{\infty} \ H_{q_i, n_i, k_i}$. The assumption that $\Phi(g)$ fixes the points in $U$ implies that it acts trivially on each factor $X_{q_i, 2, 1}$ for $i > \ell$. As each factor $H_{q_i, 2, 1}$ acts effectively on $X_{q_i, 2, 1}$ this implies that the projection of $g$ to the $i$-th factor group $H_{q_i, 2, 1}$ is the identity for $i > \ell$. This 
  implies that every entry   above the diagonal in the matrix representation of $g$ in \eqref{eq-cH} is divisible by an infinite number of distinct primes $\{q_i \mid i \geq \ell\}$, so by the Prime Factorization Theorem the matrix $g$ must be the identity. 
Alternately, observe that we have $\ds g \in \prod_{i=1}^{\ell} \ H_{q_i, 2, 1}$. This is a finite product of finite groups, which implies that $g \in \cH$ is a   torsion element. However, $\cH$ is torsion-free, hence $g$ must be the identity.  Thus, the action of $\cH$ on $X_{\infty}$ must be topologically free.

Finally, the above construction allows the choice of any infinite subset $\pi_f$  of  distinct primes, and there are an uncountable number such choices which are distinct up to asymptotic equivalence. Thus, by Theorem~\ref{thm-returnequivspectra} there are an uncountable number of topologically-free, wild nilpotent Cantor actions which are distinct up to return equivalence.
This completes the proof of Theorem~\ref{thm-lqaWILD}.

   {\begin{remark}\label{rmk-generalization}
   {\rm
   The constructions in Examples~\ref{ex-stable} and \ref{ex-wild} can be generalized to the integer upper triangular matrices in all dimensions, where there is much more freedom in the choice of the subgroups $H_{q_i, n_i,k_i}$. The above calculations become correspondingly more tedious, but yield analogous results. It seems reasonable to expect that similar constructions can be made for any finitely-generated torsion-free nilpotent (non-abelian) group $\G$.   That is, that there are group chains in $\G$ which yield wild nilpotent Cantor   actions. Note that in the work \cite{HLvL2020} with van Limbeek, the authors showed that if $\G$ is a finitely-generated nilpotent group which admits a proper self-embedding (said to be non-co-Hopfian, or renormalizable), then the iterated images of this self-embedding define a group chain for which the associated profinite action is quasi-analytic. Thus, wild Cantor actions are in a sense the furthest extreme from the actions associated to renormalizable groups.
   }
   \end{remark}
  
 \vfill
 \eject
 
%%%%%%%%%%%%%%%%%%%%%%%%%%%%%%%%%%%%%%%%%%%%%%%%%%%%%%%


\begin{thebibliography}{10}

 \bibitem{AartsFokkink1991}
{J.M.~Aarts and R.J.~Fokkink},
\newblock {\it The classification of solenoids},
\newblock {\bf Proc. Amer. Math. Soc.}, 111 :1161-1163, 1991.


  \bibitem{ALC2009}
{J.~{\'A}lvarez L{\'o}pez and A.~Candel},
\newblock {\it Equicontinuous foliated spaces},
\newblock {\bf Math. Z.}, 263:725--774, 2009.

\bibitem{ALM2016}
{J.~{\'A}lvarez L{\'o}pez and M.~Moreira Galicia},
\newblock {\it Topological {M}olino's theory},
\newblock {\bf Pacific. J. Math.}, 280:257--314, 2016.

 
\bibitem{ALBLLN2020}
{J.~{\'A}lvarez L{\'o}pez, R.~ Barral, O.~Lukina and H.~Nozawa},
\newblock {\it Wild {C}antor actions},
\newblock  {\bf J. Math. Soc. Japan},  to appear, 2021;  {arXiv:2010.00498}.



\bibitem{Auslander1988}
{J.~Auslander},
\newblock {\bf Minimal flows and their extensions},
\newblock {North-Holland Mathematics Studies}, Vol. 153, {North-Holland Publishing Co., Amsterdam}, 1988.

 
\bibitem{Baer1956}
{R.~Baer},
\newblock {\it Noethersche {G}ruppen},
\newblock {\bf Math. Z.}, 66:269--288, 1956.      

\bibitem{BartholdiGrigorchuk2000}
{L.~Bartholdi  and R.I.~Grigorchuk},
\newblock {\it On the spectrum of {H}ecke type operators related to some fractal groups},
\newblock {\bf Tr. Mat. Inst. Steklova}, 231:5--45, 2000.      

\bibitem{BartholdiGrigorchuk2002}
{L.~Bartholdi  and R.I.~Grigorchuk},
\newblock {\it On parabolic subgroups and {H}ecke algebras of some fractal groups},
\newblock {\bf Serdica Math. J.}, 28:47--90, 2002.      

  
       
\bibitem{BGS2012}
{L.~Bartholdi, R.~Grigorchuk and Z.~\v{S}uni\'{k}}, 
\newblock {\it Branch groups}, 
\newblock {\bf {Handbook of algebra, {V}ol. 3}}, {Elsevier/North-Holland, Amsterdam}, 2012, pages 989--1112.

      
\bibitem{Bing1960}
{R.H.~Bing},
\newblock {\it A simple closed curve is the only homogeneous bounded plane continuum that contains an arc},
\newblock {\bf Canad. J. Math.}, 12:209--230, 1960.

 

  \bibitem{Boyle1983}
{M.~Boyle},
\newblock {\bf Topological orbit equivalence and factor maps in symbolic dynamics},
\newblock {Ph.D. Thesis}, University of Washington, 1983.     

 \bibitem{BoyleTomiyama1998}
{M.~Boyle and J.~Tomiyama},
\newblock {\it Bounded topological orbit equivalence and {$C^*$}-algebras},
\newblock {\bf J. Math. Soc. Japan}, 50:317--329, 1998.     

\bibitem{CandelConlon2000}
{A.~Candel and L.~Conlon},
\newblock {\bf Foliations I},
\newblock Amer. Math. Soc., Providence, RI, 2000.


\bibitem{ClarkHurder2013}
{A.~Clark and S.~Hurder},
\newblock {\it Homogeneous matchbox manifolds},
\newblock {\bf Trans. A.M.S.}, 365:3151--3191, 2013.

\bibitem{CHL2019}
{A.~Clark, S.~Hurder and O.~Lukina},
\newblock {\it Classifying matchbox manifolds},
\newblock {\bf Geom \& Top},  23:1-38, 2019;  {arXiv:1311.0226}.

   \bibitem{CortezPetite2008}
{M.-I.~Cortez and S.~Petite},
\newblock {\it $G$-odometers and their almost one-to-one extensions},
\newblock {\bf J. London Math. Soc.}, 78(2):1--20, 2008. 


  \bibitem{CortezMedynets2016}
{M.I.~Cortez and K.~Medynets},
\newblock {\it Orbit equivalence rigidity of equicontinuous systems},
\newblock {\bf J. Lond. Math. Soc. (2)}, 94:545--556, 2016.

 
  \bibitem{DL2016}
{A.~Danilenko and M.~Lema\'{n}czyk},
\newblock {\it Odometer actions of the {H}eisenberg group},
\newblock {\bf J. Anal. Math.}, 128:107--157, 2016.

        
   \bibitem{DdSMS1999}
{J.D.~Dixon, M.~du~Sautoy, A.~Mann, and D.~Segal},
\newblock{\bf Analytic pro-{$p$} groups}, 
\newblock {Cambridge Studies in Advanced Mathematics}, Vol. 61, 2nd Edition, {Cambridge University Press, Cambridge}, 1999.
       
   \bibitem{dSSS2000}
{M.~du~Sautoy, D.~Segal and A.~Shalev},
\newblock{\bf New horizons in pro-{$p$} groups}, 
\newblock {Progress in Mathematics}, Vol. 184, {Birkh\"{a}user Boston, Inc., Boston, MA}, 2000.
   
     
  \bibitem{DudkoGrigorchuk2017}
{A.~Dudko and R.~Grigorchuk},
\newblock {\it On irreducibility and disjointness of {K}oopman and quasi-regular representations of weakly branch groups},
\newblock {\bf Modern theory of dynamical systems}, 
\newblock {Contemp. Math. Vol.692}, 94:51--66, 2017.

 

 \bibitem{DHL2016a}
{J.~Dyer, S.~Hurder and O.~Lukina},
\newblock {\it The discriminant invariant of Cantor group actions},
\newblock {\bf Topology Appl.}, 208: 64--92, 2016.

\bibitem{DHL2016b}
{J.~Dyer, S.~Hurder and O.~Lukina},
\newblock {\it Growth and homogeneity of matchbox manifolds},
\newblock {\bf Indag. Math.}, 28:145--169, 2017.


\bibitem{DHL2016c}
{J.~Dyer, S.~Hurder and O.~Lukina},
\newblock {\it Molino theory for matchbox manifolds},
\newblock {\bf Pacific J. Math.}, 289:91-151, 2017.
 
\bibitem{FO2002}
{R.~Fokkink and L.~Oversteegen},
\newblock {\it Homogeneous weak solenoids},
\newblock {\bf Trans. Amer. Math. Soc.}, 354(9):3743--3755, 2002.
  
\bibitem{GPS2019}
{T.~Giordano, I.~Putman and C.~Skau},
\newblock {\it {$\mathbb Z^d$}-odometers and cohomology},
\newblock {\bf Groups Geom. Dyn.}, 13:909--938, 2019.  
 
 
\bibitem{Grigorchuk2011}
{R.I.~Grigorchuk},
\newblock {\it Some problems of the dynamics of group actions on rooted trees},
{\bf Proc. Steklov Institute of Math.}, 273: 64--175, 2011.

  
 \bibitem{GL2019}
{M.~Gr\"oger and O.~Lukina},
\newblock {\it Measures and regularity of group Cantor actions}, 
\newblock {\bf  Discrete Contin. Dynam. Sys.}, to appear;  {arXiv:1911.00680}. 
 

      
\bibitem{Haefliger1984}
{A.~Haefliger},
\newblock {\it Groupo{\"\i}des d'holonomie et classifiants},
\newblock in {\bf Transversal structure of foliations (Toulouse, 1982)},
\newblock {Asterisque, 177-178, Soci\'et\'e Math\'ematique de France}, 1984:70--97.
 
 
 \bibitem{Haefliger2002a}
{A.~Haefliger},
\newblock {\it Foliations and compactly generated pseudogroups}, 
\newblock in {\bf Foliations: geometry and dynamics (Warsaw, 2000)},
\newblock {World Sci. Publ., River Edge, NJ},   2002:275--295.

  
\bibitem{HL2018a}
{S.~Hurder and O.~Lukina},
\newblock {\it Wild solenoids},
\newblock {\bf Transactions A.M.S.}, 371:4493-4533, 2019; {arXiv:1702.03032}.
 
\bibitem{HL2018b}
{S.~Hurder and O.~Lukina},
\newblock {\it Orbit equivalence and   classification of weak solenoids},
\newblock {\bf Indiana Univ. Math. J.},  Vol. 69:2339--2363, 2020; {arXiv:1803.02098}.
 
  	
 \bibitem{HL2019a}
{S.~Hurder and O.~Lukina},
\newblock {\it Limit group invariants for non-free Cantor actions},
\newblock {\bf Ergodic Theory Dynam. Systems}, to appear; {arXiv:1904.11072}.
   
\bibitem{HL2020a}
{S.~Hurder and O.~Lukina},
\newblock {\it  Nilpotent   Cantor actions},
\newblock {submitted}, 2020; {arXiv:1905.07740}.
  
 
  
\bibitem{HLvL2020}
{S.~Hurder, O.~Lukina and W.~van Limbeek},
\newblock {\it Cantor dynamics of renormalizable groups},
\newblock {\it submitted}; {arXiv:2002.01565}.
   
   
      
\bibitem{Kionke2019}
{S.~Kionke},
\newblock {\it Groups acting on rooted trees and their representations on the boundary},
\newblock {\bf J. Algebra}, 528:260--284, 2019.
    
     
  \bibitem{Li2018}
{X.~Li},
\newblock {\it Continuous orbit equivalence rigidity},
\newblock {\bf Ergodic Theory Dynam. Systems}, 38:1543--1563, 2018.

     
     
\bibitem{LSU2014}
   {S.~Lightwood, A.~ {\c{S}}ahin and I.~ Ugarcovici},
\newblock{\it The structure and spectrum of {H}eisenberg odometers},
\newblock {\bf Proc. Amer. Math. Soc.}, 142(7):{2429--2443}, 2014.

\bibitem{MS2006}
{C.C.~Moore and C.~Schochet},
\newblock {\bf Analysis on Foliated Spaces},
\newblock Math. Sci. Res. Inst. Publ. vol. 9, Second Edition,
\newblock Cambridge University Press, New York, 2006.



\bibitem{McCord1965}
{M.C.~McCord},
\newblock {\it Inverse limit sequences with covering maps},
\newblock {\bf Trans. Amer. Math. Soc.}, 114:197--209, 1965
 
 \bibitem{Nekrashevych2005}
{V.~Nekrashevych},
\newblock {\bf Self-similar groups},
\newblock {Mathematical Surveys and Monographs}, Vol. 117, 
\newblock {American Mathematical Society, Providence, RI}, 2005.  
     
 
\bibitem{Nekrashevych2016}
{V.~Nekrashevych},
\newblock {\it Palindromic subshifts and simple periodic groups of intermediate growth},
\newblock {\bf Ann. of Math. (2)}, 187:667--719, 2018.


 
 \bibitem{Renault2008}
{J.~Renault},
\newblock {\it Cartan subalgebras in {$C^*$}-algebras},
 \newblock {\bf Irish Math. Soc. Bull.}, 61:29--63, 2008.



 \bibitem{RZ2000} 
{L.~Ribes and P.~Zalesskii}, 
\newblock {\bf Profinite groups}, Springer-Verlag, Berlin 2000.


\bibitem{Sullivan2014}
{D.~Sullivan},
\newblock{Solenoidal manifolds},
 \newblock {\bf J. Singul.}, 9:203--205, 2014.


    
\bibitem{vanDantzig1930}
{D.~van Dantzig},
\newblock {\it  \"{U}ber topologisch homogene {K}ontinua},
\newblock {\bf Fund. Math.}, 15:102--125, 1930.

      
\bibitem{Wilson1998}
{J.S.~Wilson},
\newblock {\bf  Profinite groups},
\newblock {London Mathematical Society Monographs. New Series}, Vol. 19,  
\newblock The Clarendon Press, Oxford University Press, New York, 1998.  
 

\bibitem{Vietoris1927}
{L.~ Vietoris},
\newblock {\it  \"{U}ber den h\"oheren {Z}usammenhang kompakter {R}\"aume und eine {K}lasse von zusammenhangstreuen {A}bbildungen},
\newblock {\bf Math. Ann.}, 97:454--472, 1927.




        
\end{thebibliography}
\end{document}